\documentclass[12pt]{amsart}

\usepackage{graphicx,amssymb,xcolor,enumerate,latexsym}
\usepackage{hyperref}
\hypersetup{
bookmarks=true,         
pdffitwindow=false,     
pdfstartview={FitH},    
colorlinks=true,      
citecolor=red,
}

\usepackage{enumitem}

\usepackage{float}

\usepackage{mathrsfs} 
\usepackage[top=3cm, bottom=3cm, left=2.2cm, right=2.2cm]{geometry} 

\usepackage{palatino} 
\usepackage{mathpazo} 
\usepackage{amsmath} 
\usepackage{bm}

\makeatletter
\def\l@subsection{\@tocline{2}{0pt}{2.5pc}{5pc}{}}
\makeatother

\DeclareSymbolFont{largesymbol}{OMX}{yhex}{m}{n}
\DeclareMathAccent{\Widehat}{\mathord}{largesymbol}{"62}
\newcommand*\di{\mathop{}\!\mathrm{d}}

\usepackage{amsmath}
\usepackage{cleveref}

\def\B{\dot{B}_{p,q}^{3/p-1}}
\def\LB{\tilde{L}_t^\infty\dot{B}_{p,q}^{3/p-1}}
\def\LBT{\tilde{L}_t^\infty\left(0,T;\dot{B}_{p,q}^{3/p-1}\right)}

\def\LBBB{\tilde{L}_t^{\frac{2p}{p-1}}\dot{B}_{p,q}^{2/p}}
\def\LBBBT{\tilde{L}_t^{\frac{2p}{p-1}}\left(0,T;\dot{B}_{p,q}^{2/p}\right)}

\def\LBBBB{\tilde{L}_t^{\frac{p}{p-1}}\dot{B}_{p,q}^{1/p}}
\def\D{\dot{\Delta}}
\def\S{\dot{S}}
\def\R{\mathbb{R}}
\def\Z{\mathbb{Z}}
\def\A{A_4}
\def\AA{A_2}

\numberwithin{equation}{section}              
\newtheorem{theorem}{Theorem}[section]

\newtheorem{lemma}{Lemma}[section]
\newtheorem{proposition}{Proposition}[section]
\newtheorem*{proposition*}{Proposition}
\newtheorem{corollary}{Corollary}[section]
\newtheorem*{corollary*}{Corollary}
\newtheorem{definition}{Definition}[section]
\newtheorem*{definitions*}{Definitions}

\newtheorem*{acknowledgements*}{Acknowledgements}

\newtheorem*{conjecture*}{\bf Conjecture}

\newtheorem*{example*}{\bf Example}
\theoremstyle{remark}
\newtheorem{remark}{\bf Remark}[section]

\usepackage{amsthm}
\usepackage{amssymb}

\newenvironment{proof of mainresult1}[1][Proof]{\proof[\textbf{Proof of Thoerem \ref{thm:mainresult1}}]}{\endproof}

\newenvironment{proof of mainresult2}[1][Proof]{\proof[\textbf{Proof of Thoerem \ref{thm:mainresult2}}]}{\endproof}

\newenvironment{proof of analyticradius}[1][Proof]{\proof[\textbf{Proof of Corollary \ref{cor:analyticradius}}]}{\endproof}

\begin{document}

\author{Cong Wang}
\address[C. Wang]{School of Mathematical Sciences, Fudan University,  Shanghai,	200433, P.R. China.}
\email{math\_congwang@163.com}

\title[Space-time analyticity and refined analyticity radius]{Space-time analyticity and refined analyticity radius of the Navier-Stokes equations in the critical Besov spaces}

\begin{abstract}
In this paper, we establish the space-time analyticity of global solutions to the incompressible Navier-Stokes equations with small initial data in critical \emph{Besov} spaces $\B$. Time decay rates of higher order space-time joint derivatives and instantaneous lower bounds of the analyticity radius follow as straightforward consequences. The method employed combines Gevrey-class estimates with iterative derivative techniques. Furthermore,  we obtain a logarithmic improvement in the lower bound for the spatial analyticity radius of solutions for arbitrary initial data in critical \emph{Besov} spaces.
\end{abstract}

\maketitle

\section{Introduction}\label{sec:intro}
The system of incompressible Navier-Stokes equations (NSE) is given by
\begin{equation}\label{eq:NS}
\left\{
\begin{aligned}
&u_t+u\cdot\nabla u+\nabla \mathrm{p}=\mu\Delta u,\quad x\in\R^3,~~t>0,\\
&\nabla\cdot u=0,\\
&u(0,\cdot)=u_0,
\end{aligned}\right.
\end{equation}  
where $u$ is the $\mathbb{R}^3$-valued velocity field, $\mathrm{p}$ is the scalar pressure, and $\mu>0$ is the viscosity coefficient. For simplicity, we set $\mu=1$. This paper focuses on establishing the space-time analyticity and the refining of the estimate of analyticity radius of solutions to the NSE, with divergence free initial data $u_0\in\B$. A function $u:[0,T)\times\R^3\to\R^3$ is called a mild solution to the NSE, if it satisfies the integral equation:
\begin{align}\label{eq:mild}
u(t)=e^{t\Delta}u_0-\int_0^t\left[e^{(t-s)\Delta}\mathbb{P}\nabla\cdot(u\otimes u)(s)\right]\di s,
\end{align}
where $\mathbb{P}=Id-\nabla(-\Delta)^{-1}div$ is the Leray orthogonal projection and $e^{t\Delta}u_0$ is the solution to the heat equation with initial data $u_0$. Letting $u_h(t):=e^{t\Delta}u_0$, and defining the bilinear operator:
\begin{align*}
\mathcal{B}(u,v):=\int_0^t\left[e^{(t-s)\Delta}\mathbb{P}\nabla\cdot(u\otimes v)(s)\right]\di s.
\end{align*}
We rewrite $u$ as $u=u_h-\mathcal{B}(u,u)$, and further denote
\begin{align*}
\tilde{u}(t,x):=u(t,x)-u_h(t)=-\mathcal{B}(u,u).
\end{align*}
This formulation leads to the following relation:
\begin{align}\label{eq:NS1}
\tilde{u}=\mathcal{B}(\tilde{u},\tilde{u})+\mathcal{B}(\tilde{u},u_h)+\mathcal{B}(u_h,\tilde{u})+\mathcal{B}(u_h,u_h).
\end{align}
The critical spaces for the NSE, invariant under the scaling$u_{\lambda}(t,x):=\lambda u(\lambda^2t,\lambda x)$, $\mathrm{p}_\lambda(t,x):=\lambda^2\mathrm{p}(\lambda^2t,\lambda x)$, include
\begin{align*}
\dot{H}^{1/2}\hookrightarrow L^3\hookrightarrow\B\hookrightarrow BMO^{-1}\hookrightarrow\dot{B}_{\infty,\infty}^{-1},
\end{align*}
where $1<p<\infty$ and $1\leq q\leq\infty$. Kato \cite{kato1984mz} proved local well-posedness in $L^3$ and global well-posedness for small initial data. This result was later extended to $\B$ and $BMO^{-1}$ \cite{chemin1999jam,bae2012arma,koch2001advances}. It is worth noting that the NSE is ill-posed in its largest critical space $\dot{B}_{\infty,\infty}^{-1}$, which was first proved by Bourgain and Pavlovi\'{c} \cite{bourgain2008jfa}.

Analyticity of solutions to dissipative equations has important applications in fluid dynamics and optimal control theory. In fluid dynamics, the radius of space analyticity serves as a measure of the geometrically significant length scale of fluid flow \cite{grujic2001indiana} and helps in establishing upper bounds for Hausdorff lengths in the NSE \cite{kukavica1996jde}. Additionally, the analyticity of solutions contributes to the exponential convergence observed in the finite-dimensional Galerkin method applied to the Ginzburg-Landau equation \cite{doelman1993nfa0}. In optimal control theory, the analyticity of solutions can be used to obtain the observability inequalities of the dissipative equations, which is equivalent to null-controllability. For this topic, one can refer to the series of works by Wang, Zhang and their collaborators, like \cite{apraiz2014jems,escauriaza2015jmpa,wang2019jmpa}. Time analyticity of the solutions is a necessary and sufficient condition on the solvability of the backward heat equation in the class of functions with exponential growth \cite{zhang2019pams,dong2020jfa}.

In the study of the analyticity of dissipation equations, there are three particularly important methods. 
The first method is iterative derivative estimates, which is based on induction for the order of the derivative and Gronwall-type inequalities to get a coefficient growth type estimates with decay in time. This method requires a lot of computation and some technical combinatorial inequalities. Part of remarkable works can be found in \cite{germain2007imrn,dong2008arma,dong2020jfa}. In \cite{germain2007imrn}, Germain, Pavlovic, and Staffilani obtained a decay estimate in time for any spatial derivative, and space analyticity of solutions to the NSE in $BMO^{-1}$ space with small initial data. By introducing a class of spaces containing each order of spatial derivatives, Dong and Li \cite{dong2008arma} proved the space analyticity of the global solutions to subcritical dissipative quasi-geostrophic equations in its critical Lebesgue spaces. Recently, Dong and Zhang \cite{dong2020jfa} constructed a technical method to estimate the higher derivative of the Duhamel principle \eqref{eq:mild} and proved the time analyticity of the bounded solutions to the NSE. For more detail about such method, one can see \cite{wang2022jmaa,wang2023cpaa}. The second method is Gevrey estimates. This approach establishes the space analyticity of solutions by proving the boundedness of their Gevrey norm. This method originated from the work of Foias and Temam \cite{foias1989jfa}, where they studied the analyticity of periodic solutions of the NSE in both space and time with initial data $u_0\in H^1$. The result was later extended by the authors in \cite{chemin2004congress,lr2004arma}. In the critical \emph{Besov} spaces $\B$, Bae, Biswas, and Tadmor \cite{bae2012arma} obtained the space analyticity of the global solutions to the NSE with small initial data by using the Gevrey estimates method for $1\leq p<\infty$ and a special case for $p=\infty$. The third method involves extending the equations to complex spaces and using a fixed point algorithm to the complexified equations to prove space analyticity. Gruji\'{c} and Kukavica presented this idea in \cite{grujic1998jfa} for the case of initial data in $L^p$ with $3<p<\infty$. In \cite{guberovic2010dcds}, Guberovi\'{c} established space analyticity for the mild solution when $p=\infty$. Furthermore, in \cite{guberovic2010dcds}, the author also obtained the space analyticity for the global solutions in $BMO^{-1}$ with small initial data by this kind of methods. Another noteworthy result achieved through such methods is presented by Xu \cite{xu2020arma}, where the space analyticity of solutions to the NSE with initial data in $BMO$ was proved. In this paper, we propose a novel combination of Gevrey estimates for spatial variables and iterative derivative techniques for temporal variables to establish the space-time analyticity of global mild solutions with small initial data in critical $\emph{Besov}$ spaces. Additionally, we introduce a new Gevrey norm to refine the lower bound of the spatial analyticity radius.

In existing literature, most results address space and time analyticity independently. However, unlike the case of complex functions in two complex variables, where Hartogs's theorem ensures that joint analyticity follows from analyticity in each variable, the space-time analyticity of solutions to the NSE, as a function of two real variables, cannot be directly deduced by combining their separate analyticity properties. This distinction serves as a primary motivation for the present study.

The first main result of this paper is the following theorem, which asserts that global mild solutions to the NSE with small initial data in the critical \emph{Besov} spaces $\B, 1<p<\infty, 1\leq q\leq\infty$ are space-time analytic. Before stating the theorem, we introduce the functional space we are going to work with:
\begin{definition}\label{def:spaceI}
Let $1<p<\infty$, $q\leq 1\leq \infty$ and $T>0$. The space $E_{p,q}(T)$ is defined as
 \begin{align}\label{eq:Ep}
E_{p,q}(T):=\left\{f\in\mathcal{S}':~\left\|f\right\|_{E_{p,q}(T)}=\left\|f\right\|_{\LBT}+\left\|f\right\|_{\LBBBT}\right\},
\end{align}
where the definitions of spaces $\LB$ and $\LBBB$ are provided in Subsection \ref{subsec:usefullemmas}. For $T=\infty$, we denote $E_{p,q}(\infty):=E_{p,q}$. 
\end{definition}
The theorem is stated as follows:
\begin{theorem}[Space-Time Analyticity]\label{thm:mainresult1}
Let $1<p<\infty$ and $1\leq q \leq\infty$. There exists a constant $\epsilon_0>0$ such that if $u_0\in\B $ with $\left\|u_0\right\|_{\B}\leq \epsilon_0$, the mild solution $u$ to the NSE with initial data $u_0$ satisfies the following estimate
\begin{align}\label{inq:mainestimate1}
\left\|t^ne^{\sqrt{t}\Lambda}\partial_t^nu\right\|_{E_{p,q}}\leq \rho^{-1}C^nn^n,
\end{align}
for any $n\in \mathbb{N}$, where the constants $\rho,C>0$ are independent of $n$. Moreover, inequality \eqref{inq:mainestimate1} implies that $u$ is space-time analytic for any $t>0$.
\end{theorem}

A direct consequence of Theorem \ref{thm:mainresult1} is the following corollary, which provides the time decay rates of higher order space-time joint derivatives and instantaneous lower bounds on the space-time analyticity radius.

\begin{corollary}[Time Decay Rates and Analyticity Radius]\label{cor:analyticradius}
Under the conditions of Theorem \ref{thm:mainresult1}, the time decay rates of higher order derivatives of the mild solutions $u$ are given by
\begin{align*}
\left\|\partial_x^\alpha\partial_t^nu\right\|_{L^\infty}\leq C^{|\alpha|+n+2}\left(|\alpha|+n+2\right)^{|\alpha|+n+2}\left(t+1\right)t^{-\frac{|\alpha|}{2}-n-1},
\end{align*}
for any $\alpha\in\mathbb{N}^3$, $n\in\mathbb{N}$ and some constant $C>0$. Moreover, the space-time analyticity radius of $u$ satisfies
\begin{equation*}
\textrm{rad}_{S-T}(u(t))\geq\left\{
\begin{aligned}
&Ct,\quad 0<t<1;\\
&C\sqrt{t},\quad t\geq 1.
\end{aligned}\right.
\end{equation*}.
\end{corollary}

\begin{remark}\label{rmk:infinitecase}
The result of Theorem \ref{thm:mainresult1} holds for $1<p<\infty$. For $p=\infty$, the NSE  is ill-posed in $\dot{B}_{\infty, q}^{-1}$ when $q>2$, see \cite{bourgain2008jfa,germain2008jfa,yoneda2010jfa}. In \cite{bae2012arma}, the authors proved the existence of a global-in-time solution in $\dot{B}_{\infty,q}\cap\dot{B}_{3,\infty}$ with sufficiently small initial data and established that this solution is space analytic. Actually, It can be proved that such solution is also space-time analytic by combining the methods in this paper with the methods used in \cite{bae2012arma}.
\end{remark}

Another important topic concerning the analyticity of solutions to evolution equations is the refinement of the analyticity radius. Compared to the heat equation (the linear part of the NSE), $\mathcal{O}(\sqrt{t})$ seems to be the optimal lower bound for the space analyticity radius. However, Herbst and Skibsted \cite{herbst2011advances} provided the following surprising estimate for the space analyticity radius of the NSE with initial data in $H^{\gamma}$ for some $\gamma\in(1/2, 3/2)$,
\begin{align*}
\liminf_{t\to 0^+}\frac{\emph{rad}_S(u(t))}{\sqrt{t|\ln t|}}\geq\sqrt{2\gamma-1},
\end{align*}
which implies
\begin{align*}
\liminf_{t\to 0^+}\frac{\emph{rad}_S(u(t))}{\sqrt{t}}=\infty.
\end{align*}
Subsequently, this result was extended to systems of scale-invariant semi-linear parabolic equations by Chemin, Gallagher, and Zhang \cite{chemin2020mrl}. More recently, Hu and Zhang \cite{hu2022camb} established a similar result for solutions of the NSE with initial data in $L^p(\mathbb{R}^3)$ for $p\in[3,18/5)$. In this paper, we extend the result of \cite{herbst2011advances} to critical \emph{Besov} spaces, constituting our second main contribution. The functional space $E_{p,q}^\epsilon(T)$ used in this extension, is defined as follows: 
\begin{definition}\label{def:spaceII}
Let $0<\epsilon<1$, $T>0$ and $1<p<\infty$, $1\leq q\leq\infty$. The space $E_{p,q}^\epsilon(T)$ is defined as
\begin{align*}
E_{p,q}^\epsilon(T):=\left\{f\in\mathcal{S}':\left\|f\right\|_{E_{p,q}^\epsilon(T)}~\textrm{is bounded}~\right\},
\end{align*}
where
\begin{align*} \left\|f\right\|_{E_{p,q}^\epsilon(T)}=\left\|e^{-\frac{\lambda^2(t)}{4(1-\epsilon)}\frac{t}{T}}e^{\lambda(t)\frac{t}{\sqrt{T}}\Lambda}f\right\|_{\LBT}+\left\|e^{-\frac{\lambda^2(t)}{4(1-\epsilon)}\frac{t}{T}}e^{\lambda(t)\frac{t}{\sqrt{T}}\Lambda}f\right\|_{\LBBBT},
\end{align*}
with $\lambda(t)>0$ being a bounded function on $[0,T]$.
\end{definition}

Our second main result in this paper is stated as follows:
\begin{theorem}[Refined Analyticity Radius]\label{thm:mainresult2}
Let $0<\epsilon<1$, $1<p<\infty$, $1\leq q \leq\infty$ and $u_0\in\dot{B}_{p,q}^{3/p-1}$. There exists $T>0$ such that the equation \eqref{eq:NS1} has a unique solution $\tilde{u}\in E_{p,q}(T)$ with initial data $u_0$. Then, $u=u_h+\tilde{u}$ satisfies
\begin{align}\label{inq:refinedanalyticity}
\sup_{0<t\leq T}\left\|e^{-\frac{\lambda^2(t)}{4(1-\epsilon)}\frac{t}{T}}e^{\lambda(t)\frac{t}{\sqrt{T}}\Lambda}u\right\|_{\B}\leq C,
\end{align}
for some $\lambda(t)$ such that $\lim_{t\to 0}\lambda(t)=\infty$. Moreover, the space analyticity radius satisfies
\begin{align}\label{inq:refinedradius}
\liminf_{t\to 0^+}\frac{\textrm{rad}_S(u(t))}{\sqrt{t}}=\infty.
\end{align}
\end{theorem}

\begin{remark}
The space-time analyticity result in Theorem \ref{thm:mainresult1} also applies for the local solutions to the NSE with initial data of arbitrary size in $\B$. Moreover, the space-time analyticity radius also admits a refined estimate similar to \eqref{inq:refinedradius}.
\end{remark}

In this paper, the operator $\Lambda$ is a Fourier multiplier, whose symbol is given by $|\xi|_1=\sum_{i=1}^3|\xi_i|$. We emphasize that here $\Lambda\equiv\Lambda_1$ is quantified by the $\ell^1$ norm rather than by the usual $\ell^2$ norm associated with $\Lambda_2:=(-\Delta)^{1/2}$ for ease of calculation. This approach enables us to avoid cumbersome recursive estimation of higher order derivatives. In this setting, there is $\mathcal{F}\left[{e^{a\Lambda}f}\right]=e^{a|\xi|_1}\hat{f}$ for any $f\in\mathcal{S}'$ and $a\in\mathbb{R}$.

\section{Preliminaries}\label{sec:pre}
\subsection{Tool Box on Littlewood-Paley Theory}\label{subsec:littlewoodpaley}
For the convenience of the reader, we shall collect some basic facts on Littlewood-Paley theory in this subsection. Let $\mathcal{C}$ be the annulus $\left\{\xi\in\mathbb{R}^d:3/4\leq|\xi|\leq 8/3\right\}$. The smooth radial functions $\hat{\varphi}$ and $\chi$, valued in the interval $[0,1]$, belong respectively to $\mathcal{D}(\mathcal{C})$ and $\mathcal{D}(B(0,4/3))$, and are such  that
\begin{align*}
\sum_{j\in\mathbb{Z}}\hat{\varphi}(2^{-j}\xi)=1,\quad\forall\xi\in\mathbb{R}^d\backslash\{0\},
\end{align*}
and
\begin{align*}
\chi(\xi)+\sum_{j\geq 0}\hat{\varphi}(2^{-j}\xi)=1,\quad\forall\xi\in\mathbb{R}^d.
\end{align*}
Using the above two functions, we introduce the following dyadic operators:
\begin{align*}
\D_jf=\hat{\varphi}(2^{-j}D)f=2^{jd}\int_{\R^d}\varphi(2^jy)f(x-y)dy,\quad j\in\Z,
\end{align*}
\begin{align*}
\S_jf=\chi(2^{-j}D)f=2^{jd}\int_{\R^d}h(2^jy)f(x-y)dy,\quad j\in\Z,
\end{align*}
where $\varphi(x)=\mathcal{F}^{-1}[\hat{\varphi}(\xi)]$ and $h(x)=\mathcal{F}^{-1}[\chi(\xi)]$. Then, we can define the homogeneous Littlewood-Paley decomposition by
\begin{align*}
f=\sum_{j\in\Z}\D_jf,\quad\forall f\in\mathcal{S}_h',
\end{align*}
where $\mathcal{S}_h':=\{f\in\mathcal{S}': \lim_{j\to-\infty}\S_jf=0\}$. By the definition of the dyadic operators, we have
\begin{align*}
\S_jf=f-\sum_{k\geq j}\D_kf=\sum_{k\leq j-1}\D_kf.
\end{align*}

To deal with the nonlinear term in the mild solution to the NSE, we will use the following decomposition:
\begin{align}\label{inq:productdecom1}
fg=\sum_{j\in\Z}\S_jf\D_jg+\sum_{j\in\Z}\S_jg\D_jf,\quad\forall f,g\in\mathcal{S}'.
\end{align}
Then, up to finitely many terms, we have
\begin{align}\label{inq:productdecom2}
\D_j(fg)=\sum_{k\geq j-2}\D_j\left(\S_kf\D_kg\right)+\sum_{k\geq j-2}\D_j\left(\D_kf\S_kg\right),\quad\forall f,g\in\mathcal{S}'.
\end{align}
We will also use the \emph{paraproduct} given as follows:
\begin{align}\label{inq:paraproduct1}
fg=T_fg+T_gf+R(f,g)
\end{align}
with
\begin{align*}
T_fg=\sum_{i\leq j-2}\D_if\D_jg=\sum_{j\in\Z}\S_{j-1}f\D_jg\quad\textrm{and}\quad R(f,g)=\sum_{|j-j'|\leq 1}\D_jf\D_{j'}g=\sum_{j\in\mathbb{Z}}\D_jf\tilde{\Delta}_{j'}g,
\end{align*}
with 
\begin{align*}
\tilde{\Delta}_{j}g=\sum_{j'=j-1}^{j'=j+1}\D_{j'}g.
\end{align*}
Then, up to finitely many terms,
\begin{align}\label{inq:paraparoduct2}
\D_j(T_fg)=\sum_{|j'-j|\leq 4}\D_j\left(\S_{j'-1}f\D_{j'}g\right)\quad\textrm{and}\quad\D_jR(f,g)=\sum_{j'\geq j-3}\D_j\left(\D_{j'}f\tilde{\Delta}_{j'}g\right).
\end{align}

Next, we introduce homogeneous \emph{Besov} spaces and time dependent homogeneous \emph{Besov} spaces. Let $s\in\R$ and $1\leq p,q,r\leq \infty$, define
\begin{align*}
\dot{B}_{p,q}^s=\left\{f\in\mathcal{S}_h':~ \left\|f\right\|_{\dot{B}_{p,q}^s}:=\left(\sum_{j\in\Z}2^{jsq}\left\|\D_jf\right\|_{L^p}^q\right)^{\frac{1}{q}}<\infty\right\},
\end{align*}
\begin{align*}
L_t^r(0,T;\dot{B}_{p,q}^s)=\left\{f\in\mathcal{S}_h':~ \left\|f\right\|_{L_t^r(0,T;\dot{B}_{p,q}^s)}:=\left\|\left(\sum_{j\in\Z}2^{jsq}\left\|\D_jf\right\|_{L^p}^q\right)^{\frac{1}{q}}\right\|_{L_t^r(0,T)}<\infty\right\},
\end{align*}
\begin{align*}
\tilde{L}_t^r(0,T;\dot{B}_{p,q}^s)=\left\{f\in\mathcal{S}_h':~ \left\|f\right\|_{\tilde{L}_t^r(0,T;\dot{B}_{p,q}^s)}:=\left(\sum_{j\in\Z}2^{jsq}\left\|\D_jf\right\|_{L_t^r(0,T;L^p)}^q\right)^{\frac{1}{q}}<\infty\right\}.
\end{align*}

We introduce the following operators \cite{lemarie2002recent,bae2012arma}:
\begin{align*}
K_1f:=\frac{1}{2\pi}\int_0^\infty e^{ix\xi}\hat{f}(\xi)\di\xi,\qquad K_{-1}f:=\frac{1}{2\pi}\int_{-\infty}^0e^{ix\xi}\hat{f}(\xi)\di\xi,
\end{align*}
and
\begin{align*}
L_{t,1}f:=f,\qquad L_{t,-1}f:=\frac{1}{2\pi}\int_{\mathbb{R}}e^{ix\xi}e^{-2t|\xi|_1}\hat{f}(\xi)\di\xi,
\end{align*}
where $f$ is a function from $\mathbb{R}$ to $\mathbb{R}$ and $t\geq 0$ is a real number. For $\alpha=(\alpha_1,\alpha_2,\alpha_3)$, $\beta=(\beta_1,\beta_2,\beta_3)\in\{-1,1\}^3$, denote the operators
\begin{align*}
Z_{t,\alpha,\beta}=K_{\beta_1}L_{t,\alpha_1\beta_1}\otimes\cdots\otimes K_{\beta_3}L_{t,\alpha_3\beta_3}\quad\textrm{and}\quad K_\alpha=K_{\alpha_1}\otimes K_{\alpha_2}\otimes K_{\alpha_3}.
\end{align*}
The tensor product in the definition of $Z_{t,\alpha,\beta}$ and $K_\alpha$ means that the $i$-th operator in the tensor product acts on the $i$-th variable of the function $f(x_1,x_2,x_3)$. Notice that the operators $Z_{t,\alpha,\beta}$ and $K_\alpha$ commute with $\Lambda$. From the H\"ormander-Mikhlin multiplier theorem, it is easy to see that $K_\alpha$ is a bounded linear operator from $L^p(\mathbb{R}^3)$ to $L^p(\mathbb{R}^3)$ for $1<p<\infty$. Moreover, for $1<p<\infty$ and $f\in L^p(\mathbb{R}^3)$, we have
\begin{align*}
\left\|L_{t,-1}f\right\|_{L_x^p(\mathbb{R}^3)}
\leq \left\|\frac{2t}{\pi(4t^2+(\cdot)^2)}\right\|_{L_x^1(\mathbb{R})}^3\left\|f\right\|_{L_x^p(\mathbb{R}^3)}
\leq\left(\frac{4}{\pi}\right)^3 \left\|f\right\|_{L_x^p(\mathbb{R}^3)},
\end{align*}
where $\frac{2t}{\pi(4t^2+x^2)}$ is the one-dimensional poisson kernel. Hence, the operator $Z_{t,\alpha,\beta}$ is bounded from $L^p(\mathbb{R}^3)$ to $L^p(\mathbb{R}^3)$ for $1<p<\infty$, independent of $t\geq 0$. Consider
\begin{align}\label{def:operatorB}
P_t(f,g):=e^{t\Lambda}\left(e^{-t\Lambda}fe^{-t\Lambda}g\right)=\int_{\mathbb{R}^3}\int_{\mathbb{R}^3}e^{ix\cdot\left(\xi+\eta\right)}e^{t\left(|\xi+\eta|_1-|\xi|_1-|\eta|_1\right)}\hat{f}(\xi)\hat{g}(\eta)\di\xi\di\eta.
\end{align}
By splitting the integration domain and robust calculation, we have
\begin{align}\label{def:Pt}
P_t(f,g)=\sum_{(\alpha,\beta,\gamma)\in\{-1,1\}^{3\times 3}}K_{\alpha_1}\otimes K_{\alpha_2}\otimes K_{\alpha_3}\left(Z_{t,\alpha,\beta}fZ_{t,\alpha,\gamma}g\right).
\end{align}
Moreover, by the properties of $K_\alpha$ and $Z_{t,\alpha,\beta}$, we have
\begin{align}\label{inq:Pt}
\left\|P_t(f,g)\right\|_{L^p}\lesssim\left\|Z_{t,\alpha,\beta}fZ_{t,\alpha,\gamma}g\right\|_{L^p}\lesssim\left\|Z_{t,\alpha,\beta}f\right\|_{L^{p_1}}\left\|Z_{t,\alpha,\gamma}g\right\|_{L^{p_2}}
\end{align}
for any $1\leq p_1, p_2\leq\infty$ satisfy $1/p=1/p_1+1/p_2$.

\subsection{Useful Inequalities and Lemmas}\label{subsec:usefullemmas}
\begin{lemma}\label{lmm:bernsteininequality}(Bernstein's inequalities)
Let $1\leq p\leq q\leq\infty$ and $k\in\mathbb{N}$, then there exists a constant $C>0$ such that
\begin{align}\label{inq:bernstein1}
C^{-(k+1)}2^{jk}\left\|\D_jf\right\|_{L^p}\leq\sup_{|\alpha|=k}\left\|\partial_x^\alpha\D_jf\right\|_{L^p}\leq C^{k+1}2^{jk}\left\|\D_jf\right\|_{L^p},
\end{align}
and 
\begin{align}\label{inq:bernstein2}
\left\|\D_jf\right\|_{L^q}\leq C 2^{jd(1/p-1/q)}\left\|\D_j f\right\|_{L^p},
\end{align}
where $\alpha\in\mathbb{N}^3$ is a multi-index and $d\geq 1$ is the space dimension.
\end{lemma}
\begin{lemma}\label{lmm:heatkernel}(Localization of the heat kernel)
Let $1\leq p\leq\infty$. Then, there exists a constant $C>0$ such that
\begin{align}\label{inq:heat kernel}
\left\|e^{t\Delta}\D_j f\right\|_{L^p}\leq Ce^{-\frac{1}{C}t2^{2j}}\left\|\D_jf\right\|_{L^p}
\end{align}
holds for any $j\in\Z$.
\end{lemma}
\begin{lemma}\label{lmm:importantcalcultion}(\cite{dong2020jfa}, Lemma 3.3)
Let $f$ and $g$ be two smooth functions on $\R$. For any integer $n\geq 1$, we have 
\begin{equation}\label{inq:importantcalculation}
\begin{aligned}
\partial_t^n\left[t^nf(t)g(t)\right]=&\sum_{j=0}^n\binom{n}{j}\partial_t^j(t^jf(t))\partial_t^{n-j}(t^{n-j}g(t))\\
&\qquad\qquad-n\sum_{j=0}^{n-1}\binom{n-1}{j}\partial_t^j(t^jf(t))\partial_t^{n-1-j}(t^{n-1-j}g(t)).
\end{aligned}
\end{equation}
\end{lemma}

The following lemma gives a combination inequality, which was first given in \cite{kahane1969arma}, one can also find a simplified version in \cite{dong2020jfa}.
\begin{lemma}\label{lmm:combinationalinequality}(\cite{kahane1969arma}, Lemma 2.1)
There exists a universal constant $C>0$ such that
\begin{align}\label{inq:combinationinequality}
\sum_{j=0}^n\binom{n}{j}j^{j-1}(n-j)^{n-j-1}\leq Cn^{n-1}
\end{align}
holds for any $n\geq 1$.
\end{lemma}

\begin{lemma}\label{lmm:L1bounded1}(\cite{bae2012arma}, Lemma 1)
Consider the operator $E:=e^{-a\Lambda}$ for $a\geq 0$. Then $E$ is either the identity operator or is an $L^1$ kernel whose $L^1$ norm is bounded independent of $a$.
\end{lemma}

\begin{lemma}\label{lmm:Lpbounded2}(\cite{bae2012arma}, Lemma 2)
The operator $E=e^{\frac{1}{2}t\Delta+\sqrt{t}\Lambda}$ is a Fourier multiplier which maps boundedly $L^p$ to $L^p$, $1<p<\infty$, and its operator norm is uniformly bounded with respect to $t\geq 0$.
\end{lemma}

\begin{lemma}(\cite{chemin2009poincare}, Lemma 2.1)\label{lmm:fixedpointtheorem}
Let $X$ be a Banach space. Let $L$ be a continuous linear map from $X$ to $X$ and $B$ be a bilinear map from $X\times X$ to $X$. We define
\begin{align*}
\left\|L\right\|_{\mathcal{L}(X)}:=\sup_{\|x\|=1}\left\|Lx\right\|\quad\textrm{and}\quad\left\|B\right\|_{\mathcal{L}(X\times X)}:=\sup_{\|x\|=\|y\|=1}\left\|B(x,y)\right\|.
\end{align*}
If $\left\|L\right\|_{\mathcal{L}(X)}<1$, then for any $x_0$ in $X$ such that
\begin{align*}
\left\|x_0\right\|_{X}<\frac{\left(1-\left\|L\right\|_{\mathcal{L}(X)}\right)^2}{4\left\|B\right\|_{\mathcal{L}(X\times X)}},
\end{align*}
the equation
\begin{align}\label{eq:abstrctequation}
x=x_0+Lx+B(x,x)
\end{align}
has a unique solution in the ball of center $0$ and radius $\frac{1-\left\|L\right\|_{\mathcal{L}(X)}}{2\left\|B\right\|_{\mathcal{L}(X\times X)}}$.
\end{lemma}

For convenience, the constants $C>0$ appearing in the inequalities in this subsection will not be distinguished in the following of this paper. We also use $C$ to represent the upper bound of the norm of Larey orthogonal projection $\mathbb{P}$ and the upper bound of the norm of the operator $Z_{t,\alpha,\beta}$ from $L^p$ to $L^p$ for $1<p<\infty$. Additionally, we use the convention $0^s=1$ for any $s\in\mathbb{R}$. When no confusion arises, we will abbreviate $Z_{t,\alpha,\beta}$ and $Z_{t,\alpha,\gamma}$ as $Z$.

\section{Space-Time Analyticity of Solutions to the Navier-Stokes Equations}\label{sec:st-analytic}

\subsection{Existence of Solutions to the NSE in Critical \emph{Besov} Spaces}\label{subsec:existence}
In this subsection, we will state the existence result in $E_{p,q}(T)$ without proof. The existence of (global or local) solutions to the NSE in the critical \emph{Besov} spaces $\B$ has been studied in many papers, such as \cite{cannone1997rmi,chemin1999jam,bae2012arma,lemarie2002recent}.

\begin{theorem}\label{thm:existence}(Existence)
Let $1<p<\infty$, $1\leq q \leq\infty$ and $u_0\in\B$.
\begin{itemize}
\item [(i)] For any $u_0\in\B $, there exists a $T>0$ such that the NSE \eqref{eq:NS} admits a unique solution $u\in E_{p,q}(T)$.
\item [(ii)] There exists a constant $\epsilon_0$ such that for all $u_0\in\B $ with $\left\|u_0\right\|_{\B}\leq \epsilon_0$, the NSE \eqref{eq:NS} admits a global mild solution $u\in E_{p,q}$. 
\end{itemize}
\end{theorem}

In Subsection \ref{subsec:s-t analyticity}, in the first step of the proof of Proposition \ref{pro:premainresult1}, we will prove that the NSE admits a solution in a subspace of $E_{p,q}$, namely
\begin{align}\label{def:gevreyspace}
\left\{f\in E_{p,q}; e^{\sqrt{t}\Lambda}f\in E_{p,q}\right\},
\end{align}
when the initial data $u_0$ is in $\B$ and $\left\|u_0\right\|_{\B}$ is small. The method used in the proof of Proposition \ref{pro:premainresult1} also applies to Theorem \ref{thm:existence} with some modifications. In fact, the proof of Theorem \ref{thm:existence} is easier than the proof of the existence of solutions in the subspace given by \eqref{def:gevreyspace}.

\subsection{Space-Time Estimates of the Heat Kernal in \emph{Besov} Spaces}\label{subsec:heatkernel}
In this subsection, we provide several estimates on the heat kernel in \emph{Besov} spaces, which imply the space-time analyticity of the solutions to heat equations with initial data in $\B$.

\begin{lemma}\label{lmm:heatkernel1}
Let $1<p<\infty$, $1\leq q\leq\infty$, $s\in\mathbb{R}$, and $n\in\mathbb{N}$. Suppose $f\in\dot{B}_{p,q}^s$. We have $t^ne^{\sqrt{t}\Lambda}\left(\partial_t^ne^{t\Delta}f\right)\in\dot{B}_{p,q}^s$ and there exists a constant $C>0$, such that
\begin{align}\label{inq:heatkernel1}
t^n\left\|e^{\sqrt{t}\Lambda}\left(\partial_t^ne^{t\Delta}\D_jf\right)\right\|_{L^p}\leq C^{n+1}n^n\left\|\D_jf\right\|_{L^p}
\end{align}
holds for all $j\in\mathbb{Z}$, where $C>0$ is independent of $j$, $n$ and $t$. Moreover,
\begin{align}\label{inq:heatkernel2}
t^n\left\|e^{\sqrt{t}\Lambda}\left(\partial_t^ne^{t\Delta}f\right)\right\|_{\dot{B}_{p,q}^s}\leq C^{n+1}n^n\left\|f\right\|_{\dot{B}_{p,q}^s},
\end{align}
where $C>0$ is independent of $n$ and $t$.
\end{lemma}
\begin{proof}
By the heat equation, we have
\begin{align}\label{inq:heatequation}
e^{\sqrt{t}\Lambda}\left(\partial_t^ne^{t\Delta}f\right)=e^{\sqrt{t}\Lambda}\left(\Delta^ne^{t\Delta}f\right)=\Delta^ne^{\frac{1}{2}t\Delta}\left(e^{\frac{1}{2}t\Delta+\sqrt{t}\Lambda}f\right).
\end{align}
For any $1\leq p\leq \infty,j\in\mathbb{Z}$, we have
\begin{align*}
\Delta^ne^{t\Delta}\D_jf=\left(\Delta e^{\frac{t}{n}\Delta}\right)^n\D_jf.
\end{align*}
By Bernstein's inequality \eqref{inq:bernstein1} and the local estimate of heat equation \eqref{inq:heat kernel} for $n$ times, we have
\begin{align*}
&\left\|\Delta^ne^{t\Delta}\D_jf\right\|_{L^p}\\
\leq&\left(3C^42^{2j}e^{-\frac{1}{C}\frac{t}{n}2^{2j}}\right)^n\left\|\D_jf\right\|_{L^p}\\
\leq&\left(3C^5\right)^n\left(\frac{t}{n}\right)^{-n}\left(\frac{1}{C}\frac{t}{n}2^{2j}e^{-\frac{1}{C}\frac{t}{n}2^{2j}}\right)^n\left\|\D_jf\right\|_{L^p}\leq \left(3C^5\right)^n\left(\frac{t}{n}\right)^{-n}\left\|\D_jf\right\|_{L^p}.
\end{align*}
Apply $\D_j$ to \eqref{inq:heatequation} and take $L^p (1<p<\infty)$ norm to obtain
\begin{align*}
\left\|\D_je^{\sqrt{t}\Lambda}\partial_t^ne^{t\Delta}f\right\|_{L^p}=&\left\|\Delta^ne^{\frac{1}{2}t\Delta}\D_j\left(e^{\frac{1}{2}t\Delta+\sqrt{t}\Lambda}f\right)\right\|_{L^p}\\
\leq&\left(3C^5\right)^nn^n\left(\frac{1}{2}t\right)^{-n}\left\|e^{\frac{1}{2}t\Delta+\sqrt{t}\Lambda}\D_jf\right\|_{L^p}.
\end{align*}
According to Lemma \ref{lmm:Lpbounded2} or the H\"ormander-Mikhlin multiplier theorem, there exists a constant $C>0$ such that
\begin{align*}
\left\|e^{\frac{1}{2}t\Delta+\sqrt{t}\Lambda}\D_jf\right\|_{L^p}\leq C\left\|\D_j f\right\|_{L^p},
\end{align*}
for $1<p<\infty$. Hence, 
\begin{align*}
t^n\left\|\D_je^{\sqrt{t}\Lambda}\partial_t^ne^{t\Delta}f\right\|_{L^p}\leq C\left(6C^5\right)^nn^n\left\|\D_jf\right\|_{L^p},~\forall 1<p<\infty.
\end{align*}
By the definition of homogeneous \emph{Besov} spaces, the following inequality
\begin{align*}
t^n\left\|e^{\sqrt{t}\Lambda}\partial_t^n\left(e^{t\Delta}f\right)\right\|_{\dot{B}_{p,q}^s}\leq C\left(6C^5\right)^nn^n\left\|f\right\|_{\dot{B}_{p,q}^s}
\end{align*}
holds for any $t>0$, $s\in\mathbb{R}$, $1<p<\infty$ and $1\leq q\leq\infty$.
\end{proof}

\begin{lemma}\label{lmm:heatkernel2}
Let $1<p<\infty$, $1\leq q\leq\infty$, and $n\in\mathbb{N}$. Suppose $f\in\B$. There exist two constants $C, C_0>0$ such that for all $t\geq 0$, the following holds for any $j\in\mathbb{Z}$:
\begin{align}\label{inq:heatkernel3}
\left\|e^{\sqrt{t}\Lambda}\left[\partial_t^n\left(t^ne^{t\Delta}\D_jf\right)\right]\right\|_{L^p}\leq CC_0^nn^n\left\|\D_jf\right\|_{L^p}.
\end{align}
Moreover, for any $0<T\leq\infty$, we have
\begin{align}\label{inq:heatkernel4}
\left\|e^{\sqrt{t}\Lambda}\left[\partial_t^n\left(t^ne^{t\Delta}f\right)\right]\right\|_{E_{p,q}(T)}\leq CC_0^nn^n\left\|f\right\|_{\B}.
\end{align}
\end{lemma}
\begin{proof}
For any $j\in\mathbb{Z}$, applying $\D_j$ to $e^{\sqrt{t}\Lambda}\left[\partial_t^n\left(t^ne^{t\Delta}f\right)\right]$ and taking $L^p$ norm to obtain
\begin{align*}
\left\|\D_je^{\sqrt{t}\Lambda}\left[\partial_t^n\left(t^ne^{t\Delta}f\right)\right]\right\|_{L^p}\leq&\sum_{i=0}^n\binom{n}{i}n^it^{n-i}\left\|e^{\sqrt{t}\Lambda}\partial_t^{n-i}e^{\frac{1}{2}t\Delta}e^{\frac{1}{2}t\Delta}\D_jf\right\|_{L^p}\\
\leq& C\sum_{i=0}^n\binom{n}{i}n^it^{n-i}(2C)^{n-i}(n-i)^{n-i}t^{-(n-i)}\left\|e^{\frac{1}{2}t\Delta}\D_jf\right\|_{L^p}\\
\leq&C e^{-\frac{1}{2C}t2^{2j}}(2C)^nn^n\|\D_jf\|_{L^p}\sum_{i=0}^n\binom{n}{i},
\end{align*}
where we used \eqref{inq:heatkernel1} in Lemma \ref{lmm:heatkernel1} in the second inequality and used Lemma \ref{lmm:heatkernel} in the last one.
Because $\sum_{i=0}^n\binom{n}{i}=2^n$, we have
\begin{align*}
\left\|\D_je^{\sqrt{t}\Lambda}\left[\partial_t^n\left(t^ne^{t\Delta}f\right)\right]\right\|_{L^p}\leq C e^{-\frac{1}{2C}t2^{2j}}(4C)^nn^n\|\D_jf\|_{L^p},
\end{align*}
which implies \eqref{inq:heatkernel3}. 

Taking $L^{\frac{2p}{p-1}}$ norm with respect to the time variable of $\left\|\D_je^{\sqrt{t}\Lambda}\left[\partial_t^n\left(t^ne^{t\Delta}f\right)\right]\right\|_{L^p}$, we obtain
\begin{align*}
\left\|\D_je^{\sqrt{t}\Lambda}\left[\partial_t^n\left(t^ne^{t\Delta}f\right)\right]\right\|_{L^{\frac{2p}{p-1}}\left(0,T;L^p\right)}\leq C2^{-j\frac{p-1}{p}}(4C)^nn^n\left\|\D_jf\right\|_{L^p}.
\end{align*}
Hence, we have
\begin{equation}\label{inq:heatn1}
\begin{aligned}
\left\|e^{\sqrt{t}\Lambda}\left[\partial_t^n\left(t^ne^{t\Delta}f\right)\right]\right\|_{\LBBBT}\leq& \left(\sum_{j\in\mathbb{Z}}2^{j\frac{2}{p}q}\left\|\D_je^{\sqrt{t}\Lambda}\left[\partial_t^n\left(t^ne^{t\Delta}f\right)\right]\right\|_{L^\frac{2p}{p-1}\left(0,T;L^p\right)}^q\right)^{\frac{1}{q}}\\
\leq& C(4C)^nn^n\left(\sum_{j\in\mathbb{Z}}2^{j\left(\frac{3}{p}-1\right)q}\left\|\D_jf\right\|_{L^p}^q\right)^{\frac{1}{q}}\\
=&C(4C)^nn^n\left\|f\right\|_{\B}.
\end{aligned}
\end{equation}
Similarly, we also have
\begin{align}\label{inq:heatn2}
\left\|e^{\sqrt{t}\Lambda}\left[\partial_t^n\left(t^ne^{t\Delta}f\right)\right]\right\|_{\LBT}\leq C(4C)^nn^n\left\|f\right\|_{\B}.
\end{align}
Combining inequalities \eqref{inq:heatn1} and \eqref{inq:heatn2}, we obtain
\begin{align*}
\left\|e^{\sqrt{t}\Lambda}\partial_t^n\left(t^ne^{t\Delta}f\right)\right\|_{E_{p,q}(T)}=&\left\|e^{\sqrt{t}\Lambda}\partial_t^n\left(t^ne^{t\Delta}f\right)\right\|_{\LBBBT}+\left\|e^{\sqrt{t}\Lambda}\partial_t^n\left(t^ne^{t\Delta}f\right)\right\|_{\LBT}\\
\leq& 2C(4C)^nn^n\left\|f\right\|_{\B},
\end{align*}
for any $0<T\leq\infty$.
\end{proof}

\subsection{Space-Time Analyticity of the NSE}\label{subsec:s-t analyticity}
In this subsection, we prove Theorem \ref{thm:mainresult1}, which implies the mild solutions with small initial data in the critical \emph{Besov} spaces to the NSE are space-time analytic. Let $u$ be a mild solution to the NSE, we denote
\begin{align}\label{inq:Fn}
F_n(u(t)):=e^{\sqrt{t}\Lambda}\left[\partial_t^n\left(t^nu(t)\right)\right].
\end{align}
Before proving Theorem \ref{thm:mainresult1}, we first prove that there exist two constants $\rho, C>0$ such that the following inequality
\begin{align}\label{inq:weaktarget}
\left\|F_n(u)\right\|_{E_{p,q}}\leq \rho^{-1}C^nn^n
\end{align}
holds for any $n\geq 0$. In fact, later, we will show that inequality \eqref{inq:weaktarget} is equivalent to \eqref{inq:mainestimate1}. From the integral equation \eqref{eq:mild}, we have
\begin{align*}
F_n(u(t))=e^{\sqrt{t}\Lambda}\left[\partial_t^n(t^ne^{t\Delta}u_0)\right]-e^{\sqrt{t}\Lambda}\left[\partial_t^n\left(t^n\mathcal{B}(u,u)\right)\right].
\end{align*}
A direct calculation shows
\begin{align*}
e^{\sqrt{t}\Lambda}\left[\partial_t^n\left(t^n\mathcal{B}(u,u)\right)\right]&=e^{\sqrt{t}\Lambda}\partial_t^n\left[t^n\int_0^te^{(t-s)\Delta}\mathbb{P}\nabla\cdot(u\otimes u)(s)\di s\right]\\
&=\sum_{k=0}^n\binom{n}{k}e^{\sqrt{t}\Lambda}\left\{\partial_t^n\int_0^t(t-s)^ke^{(t-s)\Delta}\mathbb{P}\nabla\cdot\left[s^{n-k}\left(u(s)\otimes u(s)\right)\right]\di s\right\}.
\end{align*}
By changing the variable in the integral, we get
\begin{equation}\label{inq:changevariabel}
\begin{aligned}
&e^{\sqrt{t}\Lambda}\left[\partial_t^n\left(t^n\mathcal{B}(u,u)\right)\right]\\
=&\sum_{k=0}^n\binom{n}{k}e^{\sqrt{t}\Lambda}\left\{\partial_t^n\int_0^t s^ke^{s\Delta}\mathbb{P}\nabla\cdot\left[(t-s)^{n-k}\left(u(t-s)\otimes u(t-s)\right)\right]\di s\right\}\\
=&\sum_{k=0}^n\binom{n}{k}e^{\sqrt{t}\Lambda}\left\{\partial_t^k\int_0^t s^ke^{s\Delta} \mathbb{P}\nabla\cdot G_{n-k}\left(u,u,t-s\right)\di s\right\}\\
=&\sum_{k=0}^n\binom{n}{k}e^{\sqrt{t}\Lambda}\left\{\int_0^t \partial_t^k\left[(t-s)^ke^{(t-s)\Delta}\mathbb{P}\nabla\cdot G_{n-k}\left(u,u,s\right)\right]\di s\right\},
\end{aligned}
\end{equation}
where 
\begin{align}\label{inq:Gm}
G_m(v,w,t):=\partial_t^m\left[t^m\left(v(t)\otimes w(t)\right)\right]
\end{align}
for $m\geq 0$. This technique was first introduced in \cite{dong2020jfa}.

We will prove \eqref{inq:weaktarget} in Proposition \ref{pro:premainresult1}, and the proof of Theorem \ref{thm:mainresult1} will be provided at the end of this section. Before that, we will first present some important estimates in Proposition \ref{pro:preestimateI} and \ref{pro:preestimateII}.
\begin{proposition}\label{pro:preestimateI}
Let $1<p<\infty$, $1\leq q\leq\infty$ and $v, w\in\mathcal{S}'$. Then, there exists a constant $C>0$ such that
\begin{align}\label{inq:preI}
\left\|e^{\sqrt{t}\Lambda}\partial_t^n\left[t^n\mathcal{B}(v,w)\right]\right\|_{E_{p,q}}\leq C\sum_{k=0}^n\binom{n}{k}(16C_0)^kk^{k-1}\left\|e^{\sqrt{t}\Lambda}G_{n-k}(v,w,t)\right\|_{\tilde{L}_t^{\frac{p}{p-1}}\dot{B}_{p,q}^{1/p}}
\end{align}
 holds for any $n\in\mathbb{N}$, where $C_0>0$ is given in Lemma \ref{lmm:heatkernel2}.
\end{proposition}
\begin{proof}
A direct calculation shows
\begin{equation}\label{inq:chaincalculation}
\begin{aligned}
&\left\|e^{\sqrt{t}\Lambda}\int_0^t\partial_t^k\left[(t-s)^ke^{(t-s)\Delta}\mathbb{P}\nabla\cdot f(s)\right]\di s\right\|_{L^p}\\
=&\int_0^t\left\|e^{\left(\sqrt{t}-\sqrt{s}\right)\Lambda}\partial_t^k\left[(t-s)^ke^{(t-s)\Delta}e^{\sqrt{s}\Lambda}\mathbb{P}\nabla\cdot f(s)\right]\right\|_{L^p}\di s\\
\leq&C\int_0^t\left\|e^{\left(\sqrt{t}-\sqrt{s}\right)\Lambda}\left[\sum_{i=0}^k\binom{k}{i}\left(\partial_t^i(t-s)^k\right)\left(\partial_t^{k-i}e^{(t-s)\Delta}e^{\sqrt{s}\Lambda}\nabla\cdot f(s)\right)\right]\right\|_{L^p}\di s,
\end{aligned}
\end{equation}
where we used the boundedness of the Leray projection operator from $L^p$ to $L^p$ in the last inequality.
From the heat equation, for any $g:\mathbb{R}^3\to\mathbb{R}^3$, we have
\begin{align*}
\partial_t^l\left(e^{t\Delta}g\right)=\Delta^l\left(e^{t\Delta}g\right)=e^{\frac{1}{2}t\Delta}\Delta^l\left(e^{\frac{1}{4}t\Delta}e^{\frac{1}{4}t\Delta}g\right)=2^le^{\frac{1}{2}t\Delta}\partial_t^l\left(e^{\frac{1}{4}t\Delta}e^{\frac{1}{4}t\Delta}g\right),\quad\forall l\in\mathbb{N}.
\end{align*}
Take $g(s)=e^{\sqrt{s}\Lambda}\nabla\cdot f(s)$ and we get,
\begin{align}\label{inq:heatseparate}
\partial_t^{k-i}\left(e^{(t-s)\Delta}e^{\sqrt{s}\Lambda}\nabla\cdot f(s)\right)=2^{k-i}e^{\frac{1}{2}(t-s)\Delta}\partial_t^{k-i}\left(e^{\frac{1}{4}(t-s)\Delta}e^{\frac{1}{4}(t-s)\Delta}e^{\sqrt{s}\Lambda}\nabla\cdot f(s)\right).
\end{align}
Hence, substituting \eqref{inq:heatseparate} into \eqref{inq:chaincalculation} gives
\begin{equation*}
\begin{aligned}
&\left\|e^{\sqrt{t}\Lambda}\int_0^t\partial_t^k\left[(t-s)^ke^{(t-s)\Delta}\mathbb{P}\nabla\cdot f(s)\right]\di s\right\|_{L^p}\\
\leq&C\int_0^t\left\|e^{-a\Lambda}e^{\frac{1}{2}(t-s)\Delta+\sqrt{t-s}\Lambda}\right.\\
&\qquad\qquad\qquad\quad\left.\left[\sum_{i=0}^k\binom{k}{i}2^{k-i}\left(\partial_t^i(t-s)^k\right)\left(\partial_t^{k-i}e^{\frac{1}{4}(t-s)\Delta}\left(e^{\frac{1}{4}(t-s)\Delta}e^{\sqrt{s}\Lambda}\nabla\cdot f(s)\right)\right)\right]\right\|_{L^p}\di s,
\end{aligned}
\end{equation*}
where $a=\sqrt{t-s}+\sqrt{s}-\sqrt{t}\geq 0$. By Lemma \ref{lmm:L1bounded1} and \ref{lmm:Lpbounded2}, we obtain
\begin{align*}
&\left\|e^{\sqrt{t}\Lambda}\int_0^t\partial_t^k\left[(t-s)^ke^{(t-s)\Delta}\mathbb{P}\nabla\cdot f(s)\right]\di s\right\|_{L^p}\\
\leq&C\int_0^t\left\|\sum_{i=0}^k\binom{k}{i}2^{k-i}\left(\partial_t^i(t-s)^k\right)\left(\partial_t^{k-i}e^{\frac{1}{4}(t-s)\Delta}\left(e^{\frac{1}{4}(t-s)\Delta}e^{\sqrt{s}\Lambda} \nabla\cdot f(s)\right)\right)\right\|_{L^p}\di s\\
\leq&2^kC\int_0^t\left\|\partial_t^k\left[(t-s)^ke^{\frac{1}{4}(t-s)\Delta}\left(e^{\frac{1}{4}(t-s)\Delta}e^{\sqrt{s}\Lambda}\nabla\cdot f(s)\right)\right]\right\|_{L^p}\di s.
\end{align*}
Applying \eqref{inq:heatkernel3} in Lemma \ref{lmm:heatkernel2} to the above inequality gives
\begin{align*}
&\left\|\D_je^{\sqrt{t}\Lambda}\int_0^t\partial_t^k\left[(t-s)^ke^{(t-s)\Delta}\mathbb{P}\nabla\cdot f(s)\right]\di s\right\|_{L^p}\\
\leq& C(8C_0)^kk^k2^j\int_0^t\left\|\D_je^{\frac{1}{4}(t-s)\Delta}e^{\sqrt{s}\Lambda}f(s)\right\|_{L^p}\di s.
\end{align*}
By Lemma \ref{lmm:heatkernel}, we obtain
\begin{equation}\label{inq:withouttimenorm}
\begin{aligned}
&\left\|\D_je^{\sqrt{t}\Lambda}\int_0^t\partial_t^k\left[(t-s)^ke^{(t-s)\Delta}\mathbb{P}\nabla\cdot f(s)\right]\di s\right\|_{L^p}\\
\leq&C(8C_0)^kk^k2^j\int_0^t e^{-\frac{1}{4C}(t-s)2^{2j}}\left\|e^{\sqrt{s}\Lambda}\D_j f\right\|_{L^p}\di s.
\end{aligned}
\end{equation}
Taking the $L^{\frac{2p}{p-1}}$ norm of both sides of the above inequality with respect to time, and by Young's inequality to obtain
\begin{align*}
&\left\|\D_je^{\sqrt{t}\Lambda}\int_0^t\partial_t^k\left[(t-s)^ke^{(t-s)\Delta}\mathbb{P}\nabla\cdot f(s)\right]\di s\right\|_{L_t^{\frac{2p}{p-1}}L^p}\\
\leq&C(8C_0)^kk^k\left\|e^{-\frac{1}{4C}t2^{2j}}\right\|_{L_t^{\frac{2p}{p+1}}} 2^j\left\|e^{\sqrt{t}\Lambda}\D_jf\right\|_{L_t^{\frac{p}{p-1}}L^{p}}\leq C(8C_0)^kk^k2^{-j\frac{1}{p}}\left\|e^{\sqrt{t}\Lambda}\D_jf\right\|_{L_t^{\frac{p}{p-1}}L^{p}}.
\end{align*}
Using the relation $2^k\geq k$ for any $k\geq 1$, we have
\begin{align*}
&\left\|\D_je^{\sqrt{t}\Lambda}\int_0^t\partial_t^k\left[(t-s)^ke^{(t-s)\Delta}\mathbb{P}\nabla\cdot f(s)\right]\di s\right\|_{L_t^{\frac{2p}{p-1}}L^p}\leq C(16C_0)^kk^{k-1}2^{-j\frac{1}{p}}\left\|e^{\sqrt{t}\Lambda}\D_jf\right\|_{L_t^{\frac{p}{p-1}}L^p}.
\end{align*}
Similarly, taking $L^\infty$ norm to both sides of  inequality \eqref{inq:withouttimenorm} with respect to time, and by Young's inequality to obtain
\begin{align*}
&\left\|\D_je^{\sqrt{t}\Lambda}\int_0^t\partial_t^k\left[(t-s)^ke^{(t-s)\Delta}\mathbb{P}\nabla\cdot f(s)\right]\di s\right\|_{L_t^\infty L^p}\\
\leq&C(16C_0)^kk^{k-1}\left\|e^{-\frac{1}{4}t2^{2j}}\right\|_{L_t^p}2^j\left\|e^{\sqrt{t}\Lambda}\D_jf\right\|_{L_t^{\frac{p}{p-1}}L^p}\leq C(16C_0)^kk^{k-1}2^{j\left(1-\frac{2}{p}\right)}\left\|e^{\sqrt{t}\Lambda}\D_jf\right\|_{L_t^{\frac{p}{p-1}}L^p}.
\end{align*}
Let $f=G_{n-k}(v,w,t)$ for $n\geq k\geq 0$. Using \eqref{inq:changevariabel} and the definition of $E_{p,q}$, we obtain inequalities \eqref{inq:preI}.
\end{proof}

\begin{proposition}\label{pro:preestimateII}
Let $1<p<\infty$, $1\leq q\leq\infty$ and $v, w\in\mathcal{S}'$. Then, for any $m\geq 0$, we have the following estimate
\begin{equation}\label{inq:preII}
\begin{aligned}
\left\|e^{\sqrt{t}\Lambda}G_m(v,w,t)\right\|_{\tilde{L}^{\frac{p}{p-1}}\dot{B}_{p,q}^{1/p}}\leq&C\left[\sum_{j=0}^m\binom{m}{j}\left\|F_j(v)\right\|_{E_{p,q}}\left\|F_{m-j}(w)\right\|_{E_{p,q}}\right.\\
&\qquad\left.+m\sum_{j=0}^{m-1}\binom{m-1}{j}\left\|F_j(v)\right\|_{E_{p,q}}\left\|F_{m-1-j}(w)\right\|_{E_{p,q}}\right],
\end{aligned}
\end{equation}
where $F_n$ is defined in \eqref{inq:Fn}.
\end{proposition}
\begin{proof}
For $m\geq 0$, applying inequality \eqref{inq:importantcalculation} to $e^{\sqrt{t}\Lambda}\D_iG_m(v,w,t)$, we have
\begin{align*}
e^{\sqrt{t}\Lambda}\D_iG_m(v,w,t)=&e^{\sqrt{t}\Lambda}\D_i\left\{\sum_{j=0}^m\binom{m}{j}\left[\partial_t^j(t^jv(t))\right]\otimes\left[\partial_t^{m-j}(t^{m-j}w(t))\right]\right.\\
&\qquad\left.-m\sum_{j=0}^{m-1}\binom{m-1}{j}\left[\partial_t^j(t^jv(t))\right]\otimes\left[\partial_t^{m-1-j}(t^{m-1-j}w(t))\right]\right\}\\
=&\sum_{j=0}^m\binom{m}{j}e^{\sqrt{t}\Lambda}\D_i\left[\left(e^{-\sqrt{t}\Lambda}F_j(v(t))\right)\otimes\left( e^{-\sqrt{t}\Lambda}F_{m-j}(w(t))\right)\right]\\
&-m\sum_{j=0}^{m-1}\binom{m-1}{j}e^{\sqrt{t}\Lambda}\D_i\left[\left(e^{-\sqrt{t}\Lambda} F_j(v(t))\right)\otimes\left( e^{-\sqrt{t}\Lambda}F_{m-1-j}(w(t))\right)\right].
\end{align*}
For any $i\in\mathbb{N}$, taking $L_t^{\frac{p}{p-1}}L^p$ norm of $e^{\sqrt{t}\Lambda}\D_iG_m(v,w,t)$ and using the product decomposition \eqref{inq:productdecom1}, we have
\begin{align*}
&\left\|e^{\sqrt{t}\Lambda}\D_i G_m(v,w,t)\right\|_{L_t^{\frac{p}{p-1}}L^p}\\
\leq&C\sum_{j=0}^m\binom{m}{j}\left\|\sum_{k\geq i-2}P_{\sqrt{t}}\left(\S_kF_j(v),\D_kF_{m-j}(w)\right)\right\|_{L_t^{\frac{p}{p-1}}L^p}\\
&\qquad\quad+Cm\sum_{j=0}^{m-1}\binom{m-1}{j}\left\|\sum_{k\geq i-2}P_{\sqrt{t}}\left(\D_kF_j(v), \S_kF_{m-1-j}(w)\right)\right\|_{L_t^{\frac{p}{p-1}}L^p}
=:H_{m,i}^{(1)}+H_{m,i}^{(2)},
\end{align*}
where we used the uniformly boundedness of $\D_j(j\in\mathbb{Z})$ as a linear operator from $L^p$ to $L^p$. The definition of $P_{\sqrt{t}}$ was given in \eqref{def:Pt}.
By Bernstein's inequality \eqref{inq:bernstein2} and the relation between operators $\D_k$ and $\S_k$, we have
\begin{align*}
\left\|\S_kf\right\|_{L_t^{\frac{2p}{p-1}}L^\infty}\leq&\sum_{l\leq k-1}\left\|\D_lf\right\|_{L_t^{\frac{2p}{p-1}} L^\infty}\leq C\sum_{l\leq k-1}2^{l\frac{1}{p}}2^{l\frac{2}{p}}\left\|\D_lf\right\|_{L_t^{\frac{2p}{p-1}}L^p},
\end{align*}
for any $f\in\mathcal{S'}$. 
The discrete form of H\"older's inequality shows
\begin{align*}
\left\|\S_kf\right\|_{L_t^\frac{2p}{p-1}L^\infty}\leq C2^{k\frac{1}{p}}\left\|f\right\|_{\LBBB}\leq C 2^{k\frac{1}{p}}\left\|f\right\|_{E_{p,q}}.
\end{align*}
Hence, by \eqref{inq:Pt} and H\"older's inequality w.r.t. time variable, we have
\begin{align*}
H_{m,i}^{(1)}\leq C\sum_{j=0}^m\binom{m}{j}\left[\sum_{k\geq i-2}2^{k\frac{1}{p}}\left\|\D_kF_{m-j}(w)\right\|_{L_t^{\frac{2p}{p-1}}L^p}\left\|F_j(v)\right\|_{E_{p,q}}\right].
\end{align*}
Moreover, we also have
\begin{align*}
&2^{i\frac{1}{p}}H_{m,i}^{(1)}
\leq C\sum_{j=0}^m\binom{m}{j}\left\|F_j(v)\right\|_{E_{p,q}}\left[\sum_{k\geq i-2}2^{\frac{1}{p}(i-k)}2^{k\frac{2}{p}}\left\|\D_k F_{m-j}(w)\right\|_{L_t^{\frac{2p}{p-1}}L^p}\right].
\end{align*}
Using discrete form of  Young's inequality to the right hand side of the above inequality, we have the following inequality holds,
\begin{equation}\label{inq:preII-1}
\begin{aligned}
\left[\sum_{i\in\mathbb{Z}}\left(2^{i\frac{1}{p}}H_{m,i}^{(1)}\right)^q\right]^{\frac{1}{q}}
\leq&C\sum_{j=0}^m\binom{m}{j}\left\|F_j(v)\right\|_{E_{p,q}}\left\|F_{m-j}(w)\right\|_{\LBBB}\\
\leq&C\sum_{j=0}^m\binom{m}{j}\left\|F_j(v)\right\|_{E_{p,q}}\left\|F_{m-j}(w)\right\|_{E_{p,q}}.
\end{aligned}
\end{equation}	
Similarly, we can also obtain that
\begin{equation}\label{inq:preII-2}
\begin{aligned}
\left[\sum_{i\in\mathbb{Z}}\left(2^{i\frac{1}{p}}H_{m,i}^{(2)}\right)^q\right]^{\frac{1}{q}}\leq& Cm\sum_{j=0}^{m-1}\binom{m-1}{j}\left\|F_j(v)\right\|_{E_{p,q}}\left\|F_{m-1-j}(w)\right\|_{E_{p,q}}.
\end{aligned}
\end{equation}
Combining inequalities \eqref{inq:preII-1} and \eqref{inq:preII-2} to obtain \eqref{inq:preII}.
\end{proof}

Next, we prove inequality \eqref{inq:weaktarget}.

\begin{proposition}\label{pro:premainresult1}
Let $1< p<\infty$, $1\leq q \leq\infty$ and $u_0\in\B$.  There exists a constant $\epsilon_0>0$ such that if $\left\|u_0\right\|_{\B}\leq \epsilon_0$, the mild solution $u$ with inital data $u_0$ to the NSE satisfies the following estimate
\begin{align}\label{inq:weakmain}
\left\|F_n(u)\right\|_{E_{p,q}}\leq \rho^{-1}C^nn^{n-1},
\end{align}
for any $n\in \mathbb{N}$, where $\rho,C>0$ are constants independent of $n$.
\end{proposition}
\begin{proof}
To prove inequality \eqref{inq:weakmain}, we use induction to $n\in\mathbb{N}$. The following proof will be divided into three steps.
	
\textbf{Step 1:} For $n=0$, Proposition \ref{pro:preestimateI} and \ref{pro:preestimateII} tell us that
\begin{align}\label{inq:gevreyboundedness}
\left\|e^{\sqrt{t}\Lambda}\mathcal{B}(v,w)\right\|_{E_{p,q}}\leq C\left\|e^{\sqrt{t}\Lambda}G_0(v,w,t)\right\|_{\tilde{L}_t^{\frac{p}{p-1}}\dot{B}_{p,q}^{1/p}}\leq C\left\|e^{\sqrt{t}\Lambda}v\right\|_{E_{p,q}}\left\|e^{\sqrt{t}\Lambda}w\right\|_{E_{p,q}}.
\end{align} 
Let us take $X:=\left\{f\in E_{p,q}; e^{\sqrt{t}\Lambda}f\in E_{p,q}\right\}$ with norm $\left\|f\right\|_{X}:=\left\|e^{\sqrt{t}\Lambda}f\right\|_{E_{p,q}}$, and denote
\begin{equation*}
\begin{aligned}
x_0:=\mathcal{B}\left(u_h, u_h\right),\quad Lv:=\mathcal{B}\left(v,u_h\right),\quad\textrm{and}\quad B(v,w):=\mathcal{B}\left(v,w\right).
\end{aligned}
\end{equation*}
Then, according to \eqref{inq:gevreyboundedness} and \eqref{inq:heatkernel4}, there exists a constant $C_1>1$ such that 
\begin{align*}
\|L\|_{\mathcal{L}(X)}\leq C_1\epsilon_0,\quad\|x_0\|_X\leq C_1^2\epsilon_0^2\quad\textrm{and}\quad\|B\|_{\mathcal{L}(X\times X)}\leq C_1.
\end{align*}
Moreover, the constant $C_1$ also satisfies
\begin{align*}
\left\|e^{\sqrt{t}\Lambda}u_0\right\|_{E_{p,q}}\leq C_1\left\|u_0\right\|_{\B}\leq C_1\epsilon_0,
\end{align*}
by Lemma \ref{lmm:heatkernel2}.
Because $C_1>1$, if we choose 
\begin{align*}
\epsilon_0\leq\frac{1}{6C_1^2}<\min\left\{\frac{1}{C_1},\frac{1}{6C_1^{\frac{3}{2}}}\right\},
\end{align*}
 we will obtain from Lemma \ref{lmm:fixedpointtheorem} that there exists a unique solution $\tilde{u}$ to \eqref{eq:NS1} in space $X$, with $u_0\in\B$ satisfies $\left\|u_0\right\|_{\B}\leq\epsilon_0$, and the solution satisfies
\begin{align*}
\left\|e^{\sqrt{t}\Lambda}\tilde{u}\right\|_{E_{p,q}}<\frac{1}{2C_1}.
\end{align*}
Hence, we have
\begin{align}\label{inq:n=0}
\left\|e^{\sqrt{t}\Lambda}u\right\|_{E_{p,q}}\leq \left\|e^{\sqrt{t}\Lambda}u_0\right\|_{E_{p,q}}+\left\|e^{\sqrt{t}\Lambda}\tilde{u}\right\|_{E_{p,q}}
\leq C_1\epsilon_0+\frac{1}{2C_1}\leq\frac{2}{3C_1}:=\rho^{-1},
\end{align}
which also implies that the solution $u$ is analytic in spatial variables. Notice that we can choose $C_1>1$ sufficiently large for our needed in the following.

Let $N\in\mathbb{N}$ and $N\geq1$.  We assume that \eqref{inq:weakmain} holds for all $0\leq n\leq N-1$. More precisely, we assume
\begin{align}\label{inq:assumption}
\left\|F_n(u)\right\|_{E_{p,q}}\leq \rho^{-1}C^nn^{n-1},\quad\forall 0\leq n\leq N-1,
\end{align}
with $\rho>0$ was given in \eqref{inq:n=0} and $C:=16C_0$, where $C_0$ was given in Lemma \ref{lmm:heatkernel2}. 

\textbf{Step 2:} In this step, we will prove
\begin{equation}\label{inq:G}
\begin{aligned}
\left\|e^{\sqrt{t}\Lambda}G_m(u,u,t)\right\|_{\LBBBB}
\leq C\rho^{-2}(16C_0)^mm^{m-1}
&+\delta_N(m)C\left\|F_0(u)\right\|_{E_{p,q}}\left\|F_m(u)\right\|_{E_{p,q}}
\end{aligned}
\end{equation}
holds for all $0\leq m\leq N$ and $N\geq 1$, where
\begin{align*}
\delta_c(m)=\left\{
\begin{aligned}
&1,~m=c;\\
&0,~m\neq c,
\end{aligned}\right.
\end{align*}
and $G_m$ is defined in \eqref{inq:Gm}. Actually, if $m<N$, by the assumption \eqref{inq:assumption} and inequality \eqref{inq:preII-1}, we have 
\begin{equation}\label{inq:V1}
\begin{aligned}
\left[\sum_{i\in\mathbb{Z}}\left(2^{i\frac{1}{p}}H_{m,i}^{(1)}\right)^q\right]^{\frac{1}{q}}
\leq&C\sum_{j=0}^{m}\binom{m}{j}\left[\rho^{-1}(16C_0)^{j}j^{j-1}\right]\left[\rho^{-1}(16C_0)^{m-j}(m-j)^{m-j-1}\right]\\
\leq&C\rho^{-2}(16C_0)^mm^{m-1},
\end{aligned}
\end{equation}
where we used Lemma \ref{lmm:combinationalinequality} in the second inequality.
If $m=N$, we have
\begin{equation}\label{inq:V11}
\begin{aligned}
&\left[\sum_{i\in\mathbb{Z}}\left(2^{i\frac{1}{p}}H_{m,i}^{(1)}\right)^q\right]^{\frac{1}{q}}\\
\leq&C\sum_{j=1}^{m-1}\binom{m}{j}\left[\rho^{-1}(16C_0)^{j}j^{j-1}\right]\left[\rho^{-1}(16C_0)^{m-j}(m-j)^{m-j-1}\right]+C\left\|F_0(u)\right\|_{E_{p,q}}\left\|F_m(u)\right\|_{E_{p,q}}\\
\leq&C\rho^{-2}(16C_0)^mm^{m-1}+C\left\|F_0(u)\right\|_{E_{p,q}}\left\|F_m(u)\right\|_{E_{p,q}}.
\end{aligned}
\end{equation}
Similarly, for $1\leq m\leq N$, we can also obtain that
\begin{equation}\label{inq:V2}
\begin{aligned}
\left[\sum_{i\in\mathbb{Z}}\left(2^{i\frac{1}{p}}H_{m,i}^{(2)}\right)^q\right]^{\frac{1}{q}}\leq Cm\sum_{j=0}^{m-1}\binom{m-1}{j}\left\|F_j(u)\right\|_{E_{p,q}}\left\|F_{m-1-j}(u)\right\|_{E_{p,q}}
\leq C\rho^{-2}(16C_0)^mm^{m-1}.
\end{aligned}
\end{equation}
Combining inequalities \eqref{inq:V1} (or \eqref{inq:V11}) and \eqref{inq:V2} to obtain \eqref{inq:G}.

\textbf{Step 3:} In this step, we prove \eqref{inq:weakmain} holds for $n=N$ and this will complete the proof of this proposition. Because Proposition \ref{pro:preestimateI} and inequality \eqref{inq:G} in step 2, we have
\begin{equation*}
\begin{aligned}
&\left\|e^{\sqrt{t}\Lambda}\partial_t^N\left[t^N\mathcal{B}(u,u)\right]\right\|_{E_{p,q}}\\
\leq&C\sum_{k=0}^N\binom{N}{k}(16C_0)^kk^{k-1}\left\|e^{\sqrt{t}\Lambda}G_{N-k}(u,u,t)\right\|_{\LBBBB}\\
\leq&C\sum_{k=1}^{N-1}\binom{N}{k}(16C_0)^kk^{k-1}\left[C\rho^{-2}(16C_0)^{N-k}(N-k)^{N-k-1}\right]+C\left\|F_0(u)\right\|_{E_{p,q}}\left\|F_N(u)\right\|_{E_{p,q}}.
\end{aligned}
\end{equation*}
A direct calculation shows that
\begin{align*}
&\left\|e^{\sqrt{t}\Lambda}\partial_t^N\left[t^N\mathcal{B}(u,u)\right]\right\|_{E_{p,q}}\\
\leq&C\left[\rho^{-2}(16C_0)^N\sum_{k=1}^{N-1}\binom{N}{k}k^{k-1}(N-k)^{N-k-1}+\left\|F_0(u)\right\|_{E_{p,q}}\left\|F_N(u)\right\|_{E_{p,q}}\right].
\end{align*}
As a consequence of Lemma \ref{lmm:combinationalinequality}, we obtain that
\begin{align*}
\left\|e^{\sqrt{t}\Lambda}\partial_t^N\left[t^N\mathcal{B}(u,u)\right]\right\|_{E_{p,q}}
\leq
C\left[\rho^{-2}(16C_0)^NN^{N-1}+\left\|F_0(u)\right\|_{E_{p,q}}\left\|F_N(u)\right\|_{E_{p,q}}\right],
\end{align*}
which implies that there exists a constant $C_2>0$ such that
\begin{equation*}
\begin{aligned}
\left\|F_N(u)\right\|_{E_{p,q}}\leq&\left\|e^{\sqrt{t}\Lambda}u_0\right\|_{E_{p,q}}+\left\|e^{\sqrt{t}\Lambda}\partial_t^N\left[t^N\mathcal{B}(u,u)\right]\right\|_{E_{p,q}}\\
\leq& C_1\epsilon_0+C_2\rho^{-2}(16C_0)^NN^{N-1}+C_2\rho^{-1}\left\|F_N(u)\right\|_{E_{p,q}}.
\end{aligned}
\end{equation*}
Hence, 
\begin{align*}
\left\|F_N(u)\right\|_{E_{p,q}}\leq\frac{C_1\epsilon_0+C_2\rho^{-2}}{1-C_2\rho^{-1}}(16C_0)^NN^{N-1}.
\end{align*}
We choose $C_1$ and $C_2$ to satisfy $C_2=\frac{3}{8}C_1$, and hence $C_2=\frac{1}{4}\rho$. As a consequence, there is
\begin{align*}
\left\|F_N(u)\right\|_{E_{p,q}}\leq\rho^{-1}(32C_0)^NN^{N-1},
\end{align*}
which is \eqref{inq:weakmain}.
\end{proof}

Next, we will prove inequalities \eqref{inq:weaktarget} and \eqref{inq:mainestimate1} are equivalent. This will complete the proof of Theorem \ref{thm:mainresult1}.
\begin{proof of mainresult1}\textbf{(1) If \eqref{inq:weaktarget} holds for all $n\geq 0$, we prove \eqref{inq:mainestimate1} holds for all $n\geq 0$.}
Assume $n>0$. Notice that for $j=0,1,2,\cdots,n-1$, we have
\begin{equation*}
\begin{aligned}
t^j\partial_t^n\left(t^{n-j}u(t)\right)=&t^j\left[\sum_{i=0}^n\binom{n}{i}\partial_t^it\partial_t^{n-i}\left(t^{n-j-1}u(t)\right)\right]\\
=&t^j\left[ t\partial_t^n\left(t^{n-j-1}u(t)\right)+n\partial_t^{n-1}\left(t^{n-j-1}u(t)\right)\right]\\
=&t^{j+1}\left[\partial_t^n\left(t^{n-j-1}u(t)\right)\right]+nt^j\left[\partial_t^{n-1}\left(t^{n-j-1}u(t)\right)\right],
\end{aligned}
\end{equation*}
which implies
\begin{align}\label{eq:kinduc}
t^{j+1}\left[\partial_t^n\left(t^{n-j-1}u(t)\right)\right]=t^j\partial_t^n\left(t^{n-j}u(t)\right)-nt^j\left[\partial_t^{n-1}\left(t^{n-j-1}u(t)\right)\right].
\end{align}
Combining $j=0$ in \eqref{eq:kinduc} with inequality \eqref{inq:weaktarget}, we have
\begin{align*}
\left\|e^{\sqrt{t}\Lambda}t\partial_t^n\left(t^{n-1}u(t)\right)\right\|_{E_{p,q}}\leq&\left\|e^{\sqrt{t}\Lambda}\partial_t^n(t^nu(t))\right\|_{E_{p,q}}+\left\|ne^{\sqrt{t}\Lambda}\partial_t^{n-1}(t^{n-1}u(t))\right\|_{E_{p,q}}\\
\leq&\rho^{-1}C^n\left(1+\frac{1}{C}\right)n^n.
\end{align*}
Similarly, for any $0<k\leq n$, we have
\begin{align*}
\left\|e^{\sqrt{t}\Lambda}t\partial_t^k\left(t^{k-1}u(t)\right)\right\|_{E_{p,q}}
\leq \rho^{-1}C^k\left(1+\frac{1}{C}\right)k^k.
\end{align*}
Repeat the above process for $j=1,2,\cdots,n-1$ and we obtain
\begin{align*}
\left\|t^ne^{\sqrt{t}\Lambda}\partial_t^nu(t)\right\|_{E_{p,q}}\leq \rho^{-1}C^n\left(1+\frac{1}{C}\right)^nn^n,
\end{align*}
which implies that inequality \eqref{inq:mainestimate1} holds for $C+1$.

\textbf{(2) If \eqref{inq:mainestimate1} holds for all $n\geq 0$, we prove \eqref{inq:weaktarget} holds for all $n\geq 0$.} When $j=n-1$, from \eqref{eq:kinduc}, we obtain
\begin{align*}
\left\|t^{n-1}e^{\sqrt{t}\Lambda}\partial_t^n\left(tu(t)\right)\right\|_{E_{p,q}}\leq \left\|t^ne^{\sqrt{t}\Lambda}\partial_t^nu(t)\right\|_{E_{p,q}}+n\left\|t^{n-1}e^{\sqrt{t}\Lambda}\partial_t^{n-1}u(t)\right\|_{E_{p,q}}.
\end{align*}
As a consequence of inequality \eqref{inq:mainestimate1}, we obtain
\begin{align*}
\left\|t^{n-1}e^{\sqrt{t}\Lambda}\partial_t^n\left(tu(t)\right)\right\|_{E_{p,q}}\leq\rho^{-1}C^n\left(1+\frac{1}{C}\right)n^n.
\end{align*}
Similarly, for any $0<k\leq n$, we have
\begin{align*}
\left\|t^{k-1}e^{\sqrt{t}\Lambda}\partial_t^k\left(tu(t)\right)\right\|_{E_{p,q}}
\leq \rho^{-1}C^k\left(1+\frac{1}{C}\right)k^k.
\end{align*}
Repeat the above process for $j=n-2,n-3,\cdots,0$ and we obtain
\begin{align*}
\left\|e^{\sqrt{t}\Lambda}\partial_t^n\left(t^nu(t)\right)\right\|_{E_{p,q}}\leq \rho^{-1}C^n\left(1+\frac{1}{C}\right)^nn^n,
\end{align*}
which implies that inequality \eqref{inq:weaktarget} holds for $C+1$.
\end{proof of mainresult1}

\begin{proof of analyticradius}
By the Littlewood-Paley decomposition, we have
\begin{align*}
\left\|f\right\|_{L^\infty}\leq&\sum_{j\in\mathbb{Z}}\left\|\D_jf\right\|_{L^\infty}\leq C\sum_{j\in\mathbb{Z}}2^{j\frac{3}{p}}\left\|\D_jf\right\|_{L^p}=C\sum_{j\geq 0}2^{j\frac{3}{p}}\left\|\D_j f\right\|_{L^p}+C\sum_{j<0}2^{j\frac{3}{p}}\left\|\D_j f\right\|_{L^p}.
\end{align*}
For the first term of the right hand side of the above inequality, we use H\"older's inequality to obtain
\begin{align*}
\sum_{j\geq 0}2^{j\frac{3}{p}}\left\|\D_j f\right\|_{L^p}\leq&C\sum_{j\geq 0}2^{-j}\left[2^{j\left(\frac{3}{p}-1\right)}\sup_{|\gamma|=2}\left\|\partial_x^\gamma\D_j f\right\|_{L^p}\right]\\
\leq&C\left[\sum_{j\geq 0}\left(2^{-j}\right)^{\frac{q}{q-1}}\right]^{1-\frac{1}{q}}\sup_{|\gamma|=2}\left[\sum_{j\geq 0}2^{j\left(\frac{3}{p}-1\right)q}\left\|\partial_x^\gamma\D_j f\right\|_{L^p}^q\right]^\frac{1}{q}\\
\leq&C\sup_{|\gamma|=2}\left\|\partial_x^\gamma f\right\|_{\B}.
\end{align*}
Similarly, it is easy to proof that
\begin{align*}
\sum_{j<0}2^{j\frac{3}{p}}\left\|\D_j f\right\|_{L^p}\leq C\left\|f\right\|_{\B}.
\end{align*}
Hence, for the solution $u$ of the NSE, the following inequality holds,
\begin{align*}
\left\|\partial_x^\alpha\partial_t^nu\right\|_{L^\infty}\leq& C\sup_{|\gamma|=2}\left\|\partial_x^{\alpha+\gamma}\partial_t^nu\right\|_{\B}+C\left\|\partial_x^\alpha\partial_t^nu\right\|_{\B}\\
\leq&CC^{|\alpha|+2}(|\alpha|+2)^{|\alpha|+2}t^{-\frac{|\alpha|+2}{2}}\left\|e^{\sqrt{t}\Lambda}\partial_t^nu(t,\cdot)\right\|_{\dot{B}_{p,q}^{3/p-1}}\\
&\qquad\qquad\qquad\qquad\qquad\qquad\qquad+C C^{|\alpha|}|\alpha|^{|\alpha|}t^{-\frac{|\alpha|}{2}}\left\|e^{\sqrt{t}\Lambda}\partial_t^nu(t,\cdot)\right\|_{\dot{B}_{p,q}^{3/p-1}}.
\end{align*}
As a consequence of theorem \ref{thm:mainresult1} and Minkowski inequality, we have
\begin{align*}
\left\|\partial_x^\alpha\partial_t^nu\right\|_{L^\infty}\leq&\frac{C}{\rho} C^{|\alpha|+n+2}\left(|\alpha|+n+2\right)^{|\alpha|+n+2}t^{-\frac{|\alpha|+2}{2}-n}+\frac{C}{\rho}C^{|\alpha|+n}\left(|\alpha|+n\right)^{|\alpha|+n}t^{-\frac{|\alpha|}{2}-n}.
\end{align*}
By Stirling's formula, we obtain the instantaneous lower bound of the space-time analyticity radius as follows:
\begin{equation*}
\emph{rad}_{S-T}(t)\geq\left\{
\begin{aligned}
&Ct,\quad 0<t<1;\\
&Ct^{\frac{1}{2}},\quad t\geq 1.
\end{aligned}\right.
\end{equation*}.
\end{proof of analyticradius}

\begin{remark}
In Corollary \ref{cor:analyticradius}, when $n=0$, we obtain the instantaneous lower bound of the space analyticity radius
\begin{align*}
\emph{rad}_S(u(t))\geq C\sqrt{t}, \quad \forall t>0.
\end{align*}
Similarly, when $\alpha=0$, we obtain the instantaneous lower bound of the time analyticity radius
\begin{align*}
\emph{rad}_T(u(t))\geq Ct,\quad \forall t>0.
\end{align*}
\end{remark}

\section{Refined Analyticity Radius}\label{sec:refinedradius}
In this section, we provide a proof of Theorem \ref{thm:mainresult2}. This theorem indicates that as $t$ approaches to $0$, the space analyticity radius of the NSE with any initial data in the critical \emph{Besov} spaces $\B$ with $1<p<\infty$ and $1\leq q\leq\infty$ grows faster than $\mathcal{O}(\sqrt{t})$. This implies that the solutions to the NSE become highly analytic in space as $t$ approaches to zero. Recall Definition \ref{def:spaceII}, we will work in the space
\begin{align*}
E_{p,q}^\epsilon(T):=\left\{f\in\mathcal{S}':\left\|f\right\|_{E_{p,q}^\epsilon(T)}~\textrm{is bounded}~\right\}
\end{align*}
with
\begin{align*} \left\|f\right\|_{E_{p,q}^\epsilon(T)}=\left\|e^{-\frac{\lambda^2(t)}{4(1-\epsilon)}\frac{t}{T}}e^{\lambda(t)\frac{t}{\sqrt{T}}\Lambda}f\right\|_{\LBT}+\left\|e^{-\frac{\lambda^2(t)}{4(1-\epsilon)}\frac{t}{T}}e^{\lambda(t)\frac{t}{\sqrt{T}}\Lambda}f\right\|_{\LBBBT},
\end{align*}
where $\lambda(t)>0$ for any $t\in[0,T]$. For convenience, we denote
\begin{align*}
A_a(\lambda(s),t):=e^{-\frac{\lambda^2(s)}{a(1-\epsilon)}\frac{t}{T}},\quad A_a(t):=e^{-\frac{\lambda^2(t)}{a(1-\epsilon)}\frac{t}{T}},
\end{align*}
\begin{align*}
 f^{\lambda(s)}(t):=e^{\lambda(s)\frac{t}{\sqrt{T}}\Lambda}u(t),\quad f^\lambda(t):=e^{\lambda(t)\frac{t}{\sqrt{T}}\Lambda}f(t),
\end{align*}
and
\begin{align*}
\lambda_T:=\max_{t\in[0,T]}\lambda(t),
\end{align*}
for any $f\in\mathcal{S}'$. As a consequence, we have
\begin{align*}
\left\|f\right\|_{E_{p,q}^\epsilon(T)}=\left\|\A(t) f^{\lambda}\right\|_{\tilde{L}_t^\infty\B}+\left\|\A(t)f^{\lambda}\right\|_{\LBBB}.
\end{align*}
In the rest of this paper, we shall always use the convention that for $a\lesssim b$, we mean that there is a uniform constant $C$, which may be different on different lines, such that $a\leq Cb$. 

Next, we will present some lemmas that are the key ingredients in the proof of Theorem \ref{thm:mainresult2}. First, as a consequence of the H\"ormander-Mikhlin multiplier theorem (see Theorem 5.2.7 in \cite{grafakosclassical}), we have the following lemma:
\begin{lemma}\label{lmm:hmmultiplier}
Let $0<\epsilon<1$, $T>0$, and $1<p<\infty$. For any $\lambda(s)>0$, we have
\begin{align*}
\left\|\A(\lambda(s),t)e^{(1-\epsilon)t\Delta}e^{\lambda(s)\frac{t}{\sqrt{T}}\Lambda}u\right\|_{L^p}\leq C\left\|u\right\|_{L^p},\quad\forall s,t>0,
\end{align*}
where $C>0$ is independent of $t$, $s$ and $p$.
\end{lemma}
One can also see Lemma 3.2 in \cite{hu2022camb} to find a complete proof for Lemma \ref{lmm:hmmultiplier}.

\begin{lemma}\label{lmm:refined1}
Let $0<\epsilon<1$ and $T>0$. If $v,w\in E_{p,q}^\epsilon(T)$, then we have
\begin{align}\label{inq:refined1}
\left\|\mathcal{B}(v,w)\right\|_{E_{p,q}^\epsilon(T)}\leq Ce^{\frac{\lambda_T^2}{4(1-\epsilon)}}\left\|v\right\|_{E_{p,q}^\epsilon(T)}\left\|w\right\|_{E_{p,q}^\epsilon(T)}.
\end{align}
\end{lemma}
\begin{proof}
By the definition of $\mathcal{B}(v,w)$, for any $1<p<\infty$, we have
\begin{align*}
&\left\|\A(t) e^{\lambda(t)\frac{t}{\sqrt{T}}\Lambda}\D_j\mathcal{B}(v,w)\right\|_{L^p}\\
\lesssim& e^{\frac{\lambda_T^2}{4(1-\epsilon)}\frac{t}{T}}2^j\int_0^t\left\| e^{(t-s)\Delta}e^{\lambda(t)\frac{t}{\sqrt{T}}\Lambda}\AA(t)\D_j\left(v\otimes w\right)\right\|_{L^p}\di s\\
\lesssim&e^{\frac{\lambda_T^2}{4(1-\epsilon)}\frac{t}{T}}2^j\int_0^t\left\|M \left[e^{\epsilon(t-s)\Delta}\AA(\lambda(t),s)e^{\lambda(t)\frac{s}{\sqrt{T}}\Lambda}\D_j\left(v\otimes w\right)(s)\right]\right\|_{L^p}\di s,
\end{align*}
where $M$ is an operator defined by $Mf=\AA(\lambda(t),t-s)e^{(1-\epsilon)(t-s)\Delta}e^{\lambda(t)\frac{t-s}{\sqrt{T}}\Lambda}f$ for any function $f\in\mathcal{S}'$. Then, Lemma \ref{lmm:hmmultiplier} implies that $\left\|Mf\right\|_{L^p}\lesssim\left\|f\right\|_{L^p}$ for any $1<p<\infty$. Hence,
\begin{align*}
&\left\|\A(t) e^{\lambda(t)\frac{t}{\sqrt{T}}\Lambda}\D_j\mathcal{B}(v,w)\right\|_{L^p}\\
\lesssim&e^{\frac{\lambda_T^2}{4(1-\epsilon)}\frac{t}{T}}2^j\int_0^t\left\| e^{\epsilon(t-s)\Delta}\AA(\lambda(t),s)e^{\lambda(t)\frac{s}{\sqrt{T}}\Lambda}\D_j\left(v\otimes w\right)(s)\right\|_{L^p}\di s\\
\lesssim&e^{\frac{\lambda_T^2}{4(1-\epsilon)}\frac{t}{T}}2^j\int_0^t e^{-\epsilon(t-s)2^{2j}}\left\|\AA(\lambda(t),s)e^{\lambda(t)\frac{s}{\sqrt{T}}\Lambda}\D_j\left(v\otimes w\right)(s)\right\|_{L^p}\di s.
\end{align*}
Next, we divide the proof into three parts according to the \emph{paraproduct} \eqref{inq:paraproduct1}.

\textbf{Part I.} 
By Young's inequality w.r.t. time $t$, we have
\begin{equation}\label{inq:repart1}
\begin{aligned}
&\left\|e^{\frac{\lambda_T^2}{4(1-\epsilon)}\frac{t}{T}}2^j\int_0^t\left\| e^{-\epsilon(t-s)2^{2j}}\AA(\lambda(t),s)e^{\lambda(t)\frac{s}{\sqrt{T}}\Lambda}\D_jT_{v}w\right\|_{L^p}\di s\right\|_{L^{\infty}_t}\\
\lesssim& e^{\frac{\lambda_T^2}{4(1-\epsilon)}}2^j\left\|\AA(t)e^{\lambda(t)\frac{t}{\sqrt{T}}\Lambda}\D_jT_{v}w\right\|_{L^{\frac{p}{p-1}}_tL^p}\left\|e^{-\epsilon t2^{2j}}\right\|_{L_t^p}\\
\lesssim&e^{\frac{\lambda_T^2}{4(1-\epsilon)}}2^{j\left(1-\frac{2}{p}\right)}\left\|\AA(t)e^{\lambda(t)\frac{t}{\sqrt{T}}\Lambda}\D_jT_{v}w\right\|_{L^{\frac{p}{p-1}}_tL^p}.
\end{aligned}
\end{equation}
A direct calculation shows that
\begin{equation}\label{inq:direct}
\begin{aligned}
&\left\|\AA(t)e^{\lambda(t)\frac{t}{\sqrt{T}}\Lambda}\D_j T_{v}w(t)\right\|_{L^p}\\
\leq&\sum_{|j'-j|\leq 4}\left\|\AA(t)e^{\lambda(t)\frac{t}{\sqrt{T}}\Lambda}\D_j\left(\S_{j'-1}v(t)\D_{j'}w(t)\right)\right\|_{L^p}\\
\lesssim&\sum_{|j'-j|\leq 4}\left\|\AA(t)\D_je^{\lambda(t)\frac{t}{\sqrt{T}}\Lambda}\left[e^{-\lambda(t)\frac{t}{\sqrt{T}}\Lambda}\left(\S_{j'-1}v^{\lambda}(t)\right)e^{-\lambda(t)\frac{t}{\sqrt{T}}\Lambda}\left(\D_{j'}w^{\lambda}(t)\right)\right]\right\|_{L^p}\\
\lesssim&\sum_{|j'-j|\leq 4}\left\|\AA(t)e^{\lambda(t)\frac{t}{\sqrt{T}}\Lambda}\left[e^{-\lambda(t)\frac{t}{\sqrt{T}}\Lambda}\left(\S_{j'-1}v^{\lambda}(t)\right)e^{-\lambda(t)\frac{t}{\sqrt{T}}\Lambda}\left(\D_{j'}w^{\lambda}(t)\right)\right]\right\|_{L^p}.
\end{aligned}
\end{equation}
Recalling the definition of $P_t(f,g)$ in \eqref{def:operatorB}, we have
\begin{equation}\label{inq:Plambda}
\begin{aligned}
&\AA(t)e^{\lambda(t)\frac{t}{\sqrt{T}}\Lambda}\left[e^{-\lambda(t)\frac{t}{\sqrt{T}}\Lambda}\left(\S_{j'-1}v^{\lambda}(t)\right)e^{-\lambda(t)\frac{t}{\sqrt{T}}\Lambda}\left(\D_{j'}w^{\lambda}(t)\right)\right]\\
=&e^{\lambda(t)\frac{t}{\sqrt{T}}\Lambda}\left[\A(t)e^{-\lambda(t)\frac{t}{\sqrt{T}}\Lambda}\left(\S_{j'-1}v^{\lambda}(t)\right)\A(t)e^{-\lambda(t)\frac{t}{\sqrt{T}}\Lambda}\left(\D_{j'}w^{\lambda}(t)\right)\right]\\
=&P_{\lambda(t)\frac{t}{\sqrt{T}}}\left(\A(t)\S_{j'-1}v^{\lambda}(t),\A(t)\D_{j'}w^{\lambda}(t)\right).
\end{aligned}
\end{equation}
As a consequence of Lemma \ref{lmm:bernsteininequality}, we have
\begin{align*}
\left\|\A(t)\S_{j'-1}v^{\lambda}(t)\right\|_{L^\infty}\leq\sum_{l\leq j'-2}\left\|\A(t)\D_l v^{\lambda}(t)\right\|_{L^\infty}
\lesssim&\sum_{l\leq j'-2}2^{l\frac{3}{p}}\left\|\A(t)\D_l v^{\lambda}(t)\right\|_{L^p}\\
=&\sum_{l\leq j'-2}2^{l\frac{1}{p}}2^{l\frac{2}{p}}\left\|\A(t)\D_lv^{\lambda}(t)\right\|_{L^p}.
\end{align*}
Hence, substituting \eqref{inq:Plambda} into \eqref{inq:direct}, and as a consequence of inequality \eqref{inq:Pt}, we obtain
\begin{equation}\label{inq:useholder}
\begin{aligned}
&\left\|\AA(t)e^{\lambda(t)\frac{t}{\sqrt{T}}\Lambda}\D_j T_{v}w(t)\right\|_{L^p}\\ \lesssim&\sum_{|j'-j|\leq 4}\left\|P_{\lambda(t)\frac{t}{\sqrt{T}}}\left(\A(t)\S_{j'-1}v^{\lambda}(t),\A(t)\D_{j'}w^{\lambda}(t)\right)\right\|_{L^p}\\
\lesssim&\sum_{|j'-j|\leq 4}\left(\sum_{l\leq j'-2}2^{l\frac{1}{p}}2^{l\frac{2}{p}}\left\|\A(t)\D_lv^{\lambda}(t)\right\|_{L^p}\right)\left\|\A(t)\D_{j'}w^{\lambda}(t)\right\|_{L^p}.
\end{aligned}
\end{equation}
Hence, substituting \eqref{inq:useholder} into \eqref{inq:repart1} and then by H\"older inequality w.r.t. time and H\"older's inequality w.r.t. $\ell^1$, we have
\begin{equation}\label{inq:timeyoung}
\begin{aligned}
&\left\|e^{\frac{\lambda_T^2}{4(1-\epsilon)}\frac{t}{T}}2^j\int_0^t\left\| e^{-\epsilon(t-s)2^{2j}}\AA(\lambda(t),s)e^{\lambda(t)\frac{s}{\sqrt{T}}\Lambda}\D_jT_{v}w(s)\right\|_{L^p}\di s\right\|_{L^{\infty}_t}\\
\lesssim&e^{\frac{\lambda_T^2}{4(1-\epsilon)}}2^{j\left(1-\frac{2}{p}\right)}\sum_{|j'-j|\leq 4}\left[\left(\sum_{l\leq j'-2}2^{l\frac{1}{p}}2^{l\frac{2}{p}}\left\|\D_l\A(t)v^{\lambda}\right\|_{L_t^{\frac{2p}{p-1}} L^p}\right)\left\|\A(t)\D_{j'}w^{\lambda}\right\|_{L_t^{\frac{2p}{p-1}}L^p}\right]\\
\lesssim&e^{\frac{\lambda_T^2}{4(1-\epsilon)}}2^{j\left(1-\frac{2}{p}\right)}\left\|\A(t)v^{\lambda}\right\|_{\LBBB}\sum_{|j'-j|\leq 4}2^{j'\frac{1}{p}}\left\|\D_j\A(t) w^{\lambda}\right\|_{L_t^{\frac{2p}{p-1}}L^p}.
\end{aligned}
\end{equation}
Multiplying \eqref{inq:timeyoung} by $2^{j\left(\frac{3}{p}-1\right)}$, we obtain
\begin{align*}
&2^{j\left(\frac{3}{p}-1\right)}\left\|e^{\frac{\lambda_T^2}{4(1-\epsilon)}\frac{t}{T}}2^j\int_0^t\left\| e^{-\epsilon(t-s)2^{2j}}\AA(\lambda(t),s)e^{\lambda(t)\frac{s}{\sqrt{T}}\Lambda}\D_jT_{v}w(s)\right\|_{L^p}\di s\right\|_{L_t^\infty}\\
\lesssim&e^{\frac{\lambda_T^2}{4(1-\epsilon)}}\left\|\A(t)v^{\lambda}\right\|_{\LBBB}\sum_{|j'-j|\leq 4}2^{\left(j-j'\right)\frac{1}{p}}2^{j'\frac{2}{p}}\left\|\D_{j'}\A(t) w^{\lambda}\right\|_{L_t^{\frac{2p}{p-1}}L^p}.
\end{align*}
Hence, using the discrete form of Young's inequality w.r.t. $\ell^q$ for the above estimation, we obtain
\begin{equation}\label{inq:part1}
\begin{aligned}
&\left\|e^{\frac{\lambda_T^2}{4(1-\epsilon)}\frac{t}{T}}2^j\int_0^t e^{-\epsilon(t-s)2^{2j}}\AA(\lambda(t),s)e^{\lambda(t)\frac{s}{\sqrt{T}}\Lambda}T_{v}w(s)\di s\right\|_{\tilde{L}_t^\infty\B}\\
\lesssim& e^{\frac{\lambda_T^2}{4(1-\epsilon)}}\left\|\A(t)v^{\lambda}\right\|_{\LBBB}\left\|\A(t)w^{\lambda}\right\|_{\LBBB}.
\end{aligned}
\end{equation}
Similarly, we also have
\begin{equation}\label{inq:part2}
\begin{aligned}
&\left\|e^{\frac{\lambda_T^2}{4(1-\epsilon)}\frac{t}{T}}2^j\int_0^t e^{-\epsilon(t-s)2^{2j}}\AA(\lambda(t),s)e^{\lambda(t)\frac{s}{\sqrt{T}}\Lambda}T_{w}v(s)\di s\right\|_{\tilde{L}_t^\infty\B}\\
\lesssim& e^{\frac{\lambda_T^2}{4(1-\epsilon)}}\left\|\A(t)v^{\lambda}\right\|_{\LBBB}\left\|\A(t)w^{\lambda}\right\|_{\LBBB}.
\end{aligned}
\end{equation}

\textbf{Part II.} Young's inequality w.r.t. time imply
\begin{align*}
&2^{j\frac{2}{p}}\left\|e^{\frac{\lambda_T^2}{4(1-\epsilon)}\frac{t}{T}}2^j\int_0^t\left\| e^{-\epsilon(t-s)2^{2j}}\AA(\lambda(t),s)e^{\lambda(t)\frac{s}{\sqrt{T}}\Lambda}\D_jT_{v}w(s)\right\|_{L^p}\di s\right\|_{L_t^{\frac{2p}{p-1}}}\\
\lesssim& e^{\frac{\lambda_T^2}{4(1-\epsilon)}}2^{j\left(\frac{2}{p}+1\right)}\left\|e^{-\epsilon t2^{2j}}\right\|_{L_t^{\frac{2p}{p+1}}}\left\|\AA(t) e^{\lambda(t)\frac{t}{\sqrt{T}}\Lambda}\D_jT_{v}w\right\|_{L_t^{\frac{p}{p-1}}L^p}\\
\lesssim&e^{\frac{\lambda_T^2}{4(1-\epsilon)}}2^{j\frac{1}{p}}\left\|\AA(t) e^{\lambda(t)\frac{t}{\sqrt{T}}\Lambda}\D_jT_{v}w\right\|_{L_t^{\frac{p}{p-1}}L^p}.
\end{align*}
Hence, substituting inequality \eqref{inq:timeyoung} into the above inequality to obtain
\begin{align*}
&2^{j\frac{2}{p}}\left\|e^{\frac{\lambda_T^2}{4(1-\epsilon)}\frac{t}{T}}2^j\int_0^t e^{-\epsilon(t-s)2^{2j}}\AA(\lambda(t),s)e^{\lambda(t)\frac{s}{\sqrt{T}}\Lambda}\D_jT_{v}w(s)\di s\right\|_{L_t^{\frac{2p}{p-1}}L^p}\\
\lesssim&e^{\frac{\lambda_T^2}{4(1-\epsilon)}}\left\|A_4(t)v^{\lambda}\right\|_{\LBBB}\sum_{j'\geq j-3}2^{\left(j-j'\right)\frac{1}{p}}2^{j'\frac{2}{p}}\left\|A_4(t)\D_{j'}w^{\lambda}\right\|_{L_t^{\frac{2p}{p-1}}L^p}.
\end{align*}
Then, by Young's inequality w.r.t. $\ell^p$ space, we have
\begin{equation}\label{inq:part3}
\begin{aligned}
&\left\|e^{\frac{\lambda_T^2}{4(1-\epsilon)}\frac{t}{T}}2^j\int_0^t e^{-\epsilon(t-s)2^{2j}}\AA(\lambda(t),s)e^{\lambda(t)\frac{s}{\sqrt{T}}\Lambda}T_{v}w(s)\di s\right\|_{\LBBB}\\
\lesssim&e^{\frac{\lambda_T^2}{4(1-\epsilon)}}\left\|A_4(t)v^{\lambda}\right\|_{\LBBB}\left\|\A(t)w^{\lambda}\right\|_{\LBBB}^2.
\end{aligned}
\end{equation}
Similarly, we have
\begin{equation}\label{inq:part4}
\begin{aligned}
&\left\|e^{\frac{\lambda_T^2}{4(1-\epsilon)}\frac{t}{T}}2^j\int_0^t e^{-\epsilon(t-s)2^{2j}}\AA(\lambda(t),s)e^{\lambda(t)\frac{s}{\sqrt{T}}\Lambda}T_{w}v(s)\di s\right\|_{\LBBB}\\
\lesssim&e^{\frac{\lambda_T^2}{4(1-\epsilon)}}\left\|\A(t)v^{\lambda}\right\|_{\LBBB}\left\|\A(t)w^{\lambda}\right\|_{\LBBB}^2.
\end{aligned}
\end{equation}

\textbf{Part III.} For the term containing $R(v,w)$, by Young's inequality w.r.t. time $t$, we have
\begin{equation}\label{inq:R}
\begin{aligned}
&\left\|e^{\frac{\lambda_T^2}{4(1-\epsilon)}\frac{t}{T}}2^j\int_0^t\left\| e^{-\epsilon(t-s)2^{2j}}\AA(\lambda(t),s) e^{\lambda(t)\frac{s}{\sqrt{T}}\Lambda}\D_jR(v,w)(s)\right\|_{L^p}\di s\right\|_{L_t^\infty}\\
\lesssim&e^{\frac{\lambda_T^2}{4(1-\epsilon)}}2^{j\left(1-\frac{2}{p}\right)}\left\| \AA(t) e^{\lambda(t)\frac{t}{\sqrt{T}}\Lambda}\D_jR(v,w)\right\|_{L_t^{\frac{p}{p-1}}L^p}.
\end{aligned}
\end{equation}
H\"older's inequality and the definition and property of $P_t(f,g)$ imply
\begin{equation}\label{inq:RR}
\begin{aligned}
&\left\|\AA(t)e^{\lambda(t)\frac{t}{\sqrt{T}}\Lambda}\D_jR(v,w)(t)\right\|_{L^p}\\
\lesssim&\sum_{j'\geq j-3}\left\|\AA(t)e^{\lambda(t)\frac{t}{\sqrt{T}}\Lambda}\D_{j'}v(t)\tilde{\Delta}_{j'}w(t)\right\|_{L^p}\\
\lesssim&\sum_{j'\geq j-3}\left[\left\|\D_{j'}\A(t)v^{\lambda}(t)\right\|_{L^p}\left(\sum_{l=j'-1}^{l=j+1}2^{l\frac{3}{p}}\left\|\D_l\A(t)w^{\lambda}(t)\right\|_{L^p}\right)\right].
\end{aligned}
\end{equation}
Substituting \eqref{inq:RR} into \eqref{inq:R}, and using H\"older's inequality w.r.t. time $t$, we obtain
\begin{align*}
&\left\|e^{\frac{\lambda_T^2}{4(1-\epsilon)}\frac{t}{T}}2^j\int_0^t\left\| e^{-\epsilon(t-s)2^{2j}}\AA(\lambda(t),s) e^{\lambda(t)\frac{s}{\sqrt{T}}\Lambda}\D_jR(v,w)(s)\right\|_{L^p}\di s\right\|_{L_t^\infty}\\
\lesssim&e^{\frac{\lambda_T^2}{4(1-\epsilon)}}2^{j\left(1-\frac{2}{p}\right)}\sum_{j'\geq j-3}\left\|\D_{j'}\A(t) v^{\lambda}\right\|_{L_t^{\frac{2p}{p-1}}L^p}\left(\sum_{l=j'-1}^{l=j'+1}2^{l\frac{1}{p}}2^{l\frac{2}{p}}\left\|\D_l\A(t) w^{\lambda}\right\|_{L_t^{\frac{2p}{p-1}} L^p}\right).
\end{align*}
Using H\"older's inequality w.r.t. $\ell^q$ to obtain
\begin{align*}
&\left\|e^{\frac{\lambda_T^2}{4(1-\epsilon)}\frac{t}{T}}2^j\int_0^t\left\| e^{-\epsilon(t-s)2^{2j}}\AA(\lambda(t),s) e^{\lambda(t)\frac{s}{\sqrt{T}}\Lambda}\D_jR(v,w)(s)\right\|_{L^p}\di s\right\|_{L_t^\infty}\\
\lesssim&e^{\frac{\lambda_T^2}{4(1-\epsilon)}}\sum_{j'\geq j-3}2^{j\left(1-\frac{2}{p}\right)}2^{j'\frac{1}{p}}\left\|\D_{j'}\A(t) v^{\lambda}\right\|_{L_t^{\frac{2p}{p-1}}L^p}\left\|\A(t) w^{\lambda}\right\|_{\LBBB}.
\end{align*}
Multiplying the above inequality by $2^{j\left(\frac{3}{p}-1\right)}$. Then, we have
\begin{align*}
&2^{j\left(\frac{3}{p}-1\right)}\left\|e^{\frac{\lambda_T^2}{4(1-\epsilon)}\frac{t}{T}}2^j\int_0^t\left\| e^{-\epsilon(t-s)2^{2j}}\AA(\lambda(t),s) e^{\lambda(t)\frac{s}{\sqrt{T}}\Lambda}\D_jR(v,w)\right\|_{L^p}\di s\right\|_{L_t^\infty}\\
\lesssim&e^{\frac{\lambda_T^2}{4(1-\epsilon)}}\sum_{j'\geq j-3}2^{\left(j-j'\right)\frac{1}{p}}2^{j'\frac{2}{p}}\left\|\D_{j'}\A(t) v^{\lambda}\right\|_{L_t^{\frac{2p}{p-1}}L^p}\left\|\A(t) w^{\lambda}\right\|_{\LBBB}.
\end{align*}
Using Young's inequality to the right hand side of the above inequality w.r.t. $\ell^q$ to obtain
\begin{equation}\label{inq:part5}
\begin{aligned}
&\left\|e^{\frac{\lambda_T^2}{4(1-\epsilon)}\frac{t}{T}}\int_0^t e^{-\epsilon(t-s)2^{2j}}\AA(\lambda(t), s)e^{\lambda(t)\frac{s}{\sqrt{T}}\Lambda}R(v,w)\di s\right\|_{\tilde{L}^\infty_t\B}\\
\lesssim& e^{\frac{\lambda_T^2}{4(1-\epsilon)}}\left\|\A(t)v^{\lambda}\right\|_{\LBBB}\left\|\A(t) w^{\lambda}\right\|_{\LBBB}.
\end{aligned}
\end{equation}
Similarly, we also have
\begin{equation}\label{inq:part6}
\begin{aligned}
&\left\|e^{\frac{\lambda_T^2}{4(1-\epsilon)}\frac{t}{T}}\int_0^t e^{-\epsilon(t-s)2^{2j}}\AA(\lambda(t),s)e^{\lambda(t)\frac{s}{\sqrt{T}}\Lambda}R(v,w)(s)\di s\right\|_{\LBBB}\\
\lesssim&e^{\frac{\lambda_T^2}{4(1-\epsilon)}}\left\|\A(t) v^{\lambda}\right\|_{\LBBB}\left\|\A(t) w^{\lambda}\right\|_{\LBBB}.
\end{aligned}
\end{equation}

Finally, combining inequalities \eqref{inq:part1}, \eqref{inq:part2} and \eqref{inq:part5}, we have
\begin{align}\label{inq:partI}
\left\|\A(t) e^{\lambda(t)\frac{t}{\sqrt{T}}\Lambda}\mathcal{B}(v,w)\right\|_{\tilde{L}_t^\infty\B}\lesssim e^{\frac{\lambda_T^2}{4(1-\epsilon)}}\left\|\A(t) v^{\lambda}\right\|_{\LBBB}\left\|\A(t) w^{\lambda}\right\|_{\LBBB}.
\end{align}
Combining inequalities \eqref{inq:part3}, \eqref{inq:part4} and \eqref{inq:part6}, we also have
\begin{align}\label{inq:partII}
\left\|\A(t)e^{\lambda(t)\frac{t}{\sqrt{T}}\Lambda}\mathcal{B}(v,w)\right\|_{\LBBB}\lesssim e^{\frac{\lambda_T^2}{4(1-\epsilon)}}\left\|\A(t) v^{\lambda}\right\|_{\LBBB}\left\|\A(t) w^{\lambda}\right\|_{\LBBB}.
\end{align}
Hence, we have 
\begin{align*}
\left\|\mathcal{B}(v,w)\right\|_{E_{p,q}^\epsilon(T)}=&\left\|\A(t) e^{\lambda(t)\frac{t}{\sqrt{T}}\Lambda}\mathcal{B}(v,w)\right\|_{\tilde{L}_t^\infty\B}+\left\|\A(t)e^{\lambda(t)\frac{t}{\sqrt{T}}\Lambda}\mathcal{B}(v,w)\right\|_{\LBBB}\\
\lesssim& e^{\frac{\lambda_T^2}{4(1-\epsilon)}}\left\|v\right\|_{E_{p,q}^\epsilon(T)}\left\|w\right\|_{E_{p,q}^\epsilon(T)}.
\end{align*}
This is inequality \eqref{inq:refined1}.
\end{proof}

\begin{lemma}\label{lmm:refined2}
Let $0<\epsilon<1$ and $T>0$. If $u_o\in\B$ for $1<p<\infty$ and $1\leq q\leq\infty$. Then there holds
\begin{align}\label{inq:refined2}
\left\|\mathcal{B}(u_h,u_h)\right\|_{E_{p,q}^\epsilon(T)}\leq Ce^{\frac{\lambda_T^2}{4(1-\epsilon)}}\left\|u_h\right\|_{\LBBB}^2.
\end{align}
\end{lemma}
\begin{proof}
Let $v=w=u_h$ in Lemma \ref{lmm:refined1}. By applying inequalities \eqref{inq:partI} and \eqref{inq:partII}, we can directly derive inequality \eqref{inq:refined2}.
\end{proof}

\begin{lemma}\label{lmm:refined3}
Let $0<\epsilon<1$ and $T>0$. If $u_o\in\dot{B}_{p,q}^{3/p-1}$ for $1<p<\infty$ and $1\leq q\leq\infty$, and $v\in E_{p,q}^\epsilon(T)$, one has
\begin{align}\label{inq:refined3}
\left\|\mathcal{B}(v,u_h)\right\|_{E_{p,q}^\epsilon(T)}\leq Ce^{\frac{\lambda_T^2}{4(1-\epsilon)}}\left\|v\right\|_{E_{p,q}^\epsilon(T)}\left\|u_h\right\|_{\tilde{L}_T^{\frac{2p}{p-1}}\dot{B}_{p,q}^{2/p}}.
\end{align}
\end{lemma}
\begin{proof}
Let $w=u_h$ in Lemma \ref{lmm:refined1}. Using inequalities \eqref{inq:partI} and \eqref{inq:partII}, and the definition of $E_{p,q}^\epsilon$, we will obtain \eqref{inq:refined3} directly.
\end{proof}

With the above lemmas, we are ready to prove our second main result, Theorem \ref{thm:mainresult2}.

\begin{proof of mainresult2}
Let us take $X:=E_{p,q}^\epsilon(T)$ and 
\begin{equation*}
\begin{aligned}
x_0:=\mathcal{B}\left(u_h, u_h\right),\quad Lv:=\mathcal{B}\left(v,u_h\right),\quad\textrm{and}\quad B(v,w):=\mathcal{B}\left(v,w\right).
\end{aligned}
\end{equation*}
Moreover, according to Lemmas \ref{lmm:refined1}, \ref{lmm:refined2}, and \ref{lmm:refined3}, there exists a constant $C>0$ such that
\begin{align*}
\left\|x_0\right\|_{X}\leq Ce^{\frac{\lambda_T^2}{4(1-\epsilon)}}\left\|u_h\right\|_{\LBBB}^2,\quad\left\|L\right\|_{\mathcal{L}(X)}\leq Ce^{\frac{\lambda_T^2}{4(1-\epsilon)}}\left\|u_h\right\|_{\LBBB},
\end{align*}
and 
\begin{align*}
\left\|B\right\|_{\mathcal{L}(X\times X)}\leq Ce^{\frac{\lambda_T^2}{4(1-\epsilon)}}.
\end{align*}
We denote
\begin{align*}
\lambda(T)=2\sqrt{(1-\epsilon)\ln{\frac{1}{\left\|u_h\right\|_{\LBBB}^{\frac{1}{2}}}}},
\end{align*}
By the Weierstrass M-test, we have
\begin{align*}
\lim_{T\to 0}\left\|u_h\right\|_{\LBBB}=0,
\end{align*}
which implies that the definition of $\lambda(T)$ is compatible when $T$ is small. As a consequence, there are
\begin{align*}
\lim_{T\to 0}\lambda(T)=\infty,
\end{align*}
and
\begin{align*}
\lim_{T\to 0}e^{\frac{\lambda_T^2}{4(1-\epsilon)}}\left\|u_h\right\|_{\LBBB}=\lim_{T\to 0}e^{\frac{\lambda^2(T)}{4(1-\epsilon)}}\left\|u_h\right\|_{\LBBB}
=\lim_{T\to 0}\left\|u_h\right\|_{\LBBB}^{\frac{1}{2}}=0.
\end{align*}
Hence, we can choose some $T>0$ so small that
\begin{align*}
Ce^{\frac{\lambda_T^2}{4(1-\epsilon)}}\left\|u_h\right\|_{\LBBB}\leq 1,
\end{align*}
and 
\begin{align*}
Ce^{\frac{\lambda_T^2}{4(1-\epsilon)}}\left\|u_h\right\|_{\LBBB}^2<\frac{\left(1-\left\|L\right\|_{\mathcal{L}(X)}\right)^2}{4\left\|B\right\|_{\mathcal{B}(X)}}.
\end{align*}
As a consequence of Lemma \ref{lmm:fixedpointtheorem}, there exist some $T>0$ and a unique solution $\tilde{u}$ of equation \eqref{eq:NS1} in $E_{p,q}^\epsilon(T)$ to equation \eqref{eq:NS1}.  As a consequence of Minkowski's inequality, the following estimate 
 \begin{align*}
\sup_{0<t\leq T} \left\|e^{-\frac{\lambda^2(t)}{4(1-\epsilon)}\frac{t}{T}}e^{\lambda(t)\frac{t}{\sqrt{T}}\Lambda}\tilde{u}\right\|_{\B}\leq C
 \end{align*}
 holds for some $C>0$. Moreover, inequality \eqref{inq:refinedanalyticity} and the refined estimate \eqref{inq:refinedradius} hold because of that it can be easy to prove $u_o\in E_{p,q}^\epsilon$.
\end{proof of mainresult2}

\quad\\
\noindent\textbf{Data availability.} No data was used for the research described in the article.\\
\\
\noindent\textbf{Conflict of interest.} The authors declare that they have no conflict of interest.

\appendix

\bibliographystyle{plain}
\bibliography{NS_Besov}

\end{document}